\documentclass[11pt]{article}
\usepackage{graphicx}
\usepackage{color}
\usepackage{amsmath}
\usepackage{amssymb}
\usepackage{amscd}
\newcommand{\R}{\mathbb{R}}
\newcommand{\inr}[1]{\bigl< #1 \bigr>}

\newcommand{\E}{\mathbb{E}}

\newcommand{\eps}{\varepsilon}

\newtheorem{Theorem}{Theorem}[section]
\newtheorem{Lemma}[Theorem]{Lemma}

\newtheorem{Definition}[Theorem]{Definition}

\newtheorem{Corollary}[Theorem]{Corollary}
\newtheorem{Remark}[Theorem]{Remark}

\numberwithin{equation}{section} 

\def \proof {\noindent {\bf Proof.}\ \ }

\def \endproof
{{\mbox{}\nolinebreak\hfill\rule{2mm}{2mm}\par\medbreak}}
\def\IND{1}

\begin{document}
\title{Empirical processes with a bounded $\psi_1$ diameter}
\author{
Shahar Mendelson\footnote{Department of Mathematics, Technion,
I.I.T, Haifa 32000, Israel, and Centre for Mathematics and its
Applications, Institute of Advanced Studies, The Australian
National University, Canberra, ACT 0200, Australia.
\newline
{\sf email: shahar@tx.technion.ac.il, \ \ shahar.mendelson@anu.edu.au
}
 }
}
\medskip
\maketitle
\begin{abstract}
We study the empirical process $\sup_{f \in F} |N^{-1}\sum_{i=1}^N
f^2(X_i)-\E f^2|$, where $F$ is a class of mean-zero functions on a
probability space $(\Omega,\mu)$ and $(X_i)_{i=1}^N$ are selected
independently according to $\mu$.

We present a sharp bound on this supremum that depends on the
$\psi_1$ diameter of the class $F$ (rather than on the $\psi_2$ one)
and on the complexity parameter $\gamma_2(F,\psi_2)$. In addition,
we present optimal bounds on the random diameters $\sup_{f \in F}
\max_{|I|=m} (\sum_{i \in I} f^2(X_i))^{1/2}$ using the same
parameters. As applications, we extend several well known results in
Asymptotic Geometric Analysis to any isotropic, log-concave ensemble
on $\R^n$.
\end{abstract}
\section{Introduction} \label{sec:intro}
In this article we study the empirical process
\begin{equation} \label{eq:zero}
\sup_{f \in F} \left| \frac{1}{N} \sum_{i=1}^N f^2(X_i) - \E f^2
\right|,
\end{equation}
where $F$ is a class of functions on the probability space
$(\Omega,\mu)$ and $(X_i)_{i=1}^N$ are independent, distributed
according to $\mu$. Properties of this process play an important
part in Asymptotic Geometric Analysis and in Nonparametric
Statistics, though even without considering the possible applications,
\eqref{eq:zero} is a natural object. Indeed, a fundamental problem
in Empirical Processes Theory is to understand the way the empirical
(random) structure of a class of functions, obtained by random
sampling, captures the original structure determined by the
underlying measure $\mu$. More accurately, one wishes to relate,
with high probability, $N^{-1}\sum_{i=1}^N \ell(f(X_i))$ to $\E
\ell(f)$, uniformly in $f \in F$, for a reasonable real valued
function $\ell$. The two most natural functions that are considered
in this context are $\ell(t)=t$, which leads to the Uniform Law of
Large Numbers, and $\ell(t)=t^2$, which is connected to properties
of the Uniform Central Limit Theorem and gives information on the
way the empirical $\ell_2$ structure of $F$ is connected to the
$L_2(\mu)$ one (see \cite{Dud-book,VW} for an extensive study of
these topics).

Despite its importance, bounds on \eqref{eq:zero} are not
satisfactory. Standard empirical processes methods allow one to
bound \eqref{eq:zero} only in rather trivial cases, in which either
the class $F$ is bounded in $L_\infty$, or if it has a well behaved
envelope function (recall that an envelope function is
$W(\omega)=\sup_{f \in F} |f(\omega)|$). In those cases it is
possible to use contraction methods and control \eqref{eq:zero}
using the linear process $\sup_{f \in F} |N^{-1}\sum_{i=1}^N f(X_i)
-\E f|$, which is a far simpler object than \eqref{eq:zero}.
However, if the function class is not uniformly bounded, or even if
it is, but with a very weak uniform bound, contraction based methods
lead to trivial estimates on \eqref{eq:zero}.

An alternative approach is to control \eqref{eq:zero} using random
parameters of $F$ that depend on the geometry of typical
coordinate projections
$$
P_\sigma F =\left\{ \left(f(X_1),...,f(X_N)\right) \ : \ f \in F
\right\},
$$
for an independent sample $\sigma=(X_1,...,X_N)$. The downside of
this approach is that the structure of $P_\sigma F$ itself is often
difficult to handle, let alone that of $P_\sigma F^2$ for $F^2=\{f^2
: f \in F\}$. Moreover, the standard way of relating the geometry of
$P_\sigma F^2$ to that of $P_\sigma F$ also involves contraction
methods, resulting in the same type of problems that have been mentioned
above.

To illustrate the difficulty, consider the following, seemingly
simple problem. Let $\Omega = \R^n$ and assume that $\mu$ is a
natural measure on $\R^n$, say the canonical gaussian measure, the
uniform measure on $\{-1,1\}^n$, or more generally, an isotopic
log-concave measure (see the definitions in Section \ref{sec:pre}).
Let $F$ be the class of linear functionals on $\R^n$ of Euclidean
norm one, that is, $F=\{\inr{x,\cdot} : x \in S^{n-1}\}$.

Note that $F$ may consist of unbounded functions on $(\Omega,\mu)$,
or, at best, of functions with an $L_\infty$ bound that grows
polynomially with the dimension $n$. It is straightforward to show
that contraction based methods lead to a very loose estimate on
\eqref{eq:zero} in such a case. To make things worse, if one
considers a typical sample $(X_i)_{i=1}^N$, the structure of the
ellipsoid $P_\sigma F$ is hard to handle (certainly if all the
information that one has on $\mu$ is that it is an isotropic,
log-concave measure). And, finally, an attempt to bound \eqref{eq:zero}
using the structure of the class $F^2=\{f^2 : f \in F\}$ directly,
without linearizing, will fail because $P_\sigma F^2$ is a rather
complicated object.

It would be highly desirable to bound \eqref{eq:zero} using a
deterministic parameter of $F$, that is, a metric invariant of $F$
that depends on $\mu$ and not on $(X_i)_{i=1}^N$, since in many
applications (the example mentioned above for one), $F$ has a simple
structure relative to a natural metric. Thus, our aim here is to
obtain bounds on \eqref{eq:zero} that depend on the deterministic
structure of $F$ as a class of functions on $(\Omega,\mu)$. All we will assume is
that $F$ consists of functions that have well behaved tails, but may
be unbounded, and the class may be without a good envelope function.

It turns out that if one wishes to bound $\E \sup_{f \in F}
\left|N^{-1} \sum_{i=1}^N f(X_i) - \E f \right|$ using a
deterministic metric structure of $F$, one has to consider metrics
that are stronger than the $L_p(\mu)$ ones (see Lemma
\ref{lemma:psi-alpha-lower-bound} and Remark
\ref{rem:psi-alpha-lower-bound} for an exact formulation of this observation). More
reasonable metrics for such a goal are the Orlicz norms
$\psi_\alpha$ for $1 \leq \alpha \leq 2$. These norms are defined
via the Young function $\exp(x^\alpha)-1$ for $\alpha \geq 1$ by
$$
\|f\|_{\psi_\alpha} = \inf \left\{c>0 : \E \exp(|f/c|^\alpha) \leq 2
\right\}.
$$

It is possible to bound the ``linear" process using a
natural complexity parameter of $F$ that originated in the theory of
Gaussian Processes. This complexity parameter is defined for any
metric space $(T,d)$ and is denoted by $\gamma_2(T,d)$ (see the book
\cite{Tal:book} and Section \ref{sec:pre} for its definition and
some of its properties). Indeed, it is standard to show that
$$
\E \sup_{f \in F} \left|\frac{1}{N} \sum_{i=1}^N f(X_i) - \E f
\right| \leq c \frac{\gamma_2(F,\psi_2)}{\sqrt{N}},
$$
where $c$ is an absolute constant, and that a similar bound holds
with high probability (see Lemma \ref{thm:generic-chaining}).
Moreover, as we will explain in Section \ref{sec:diam-optimality},
it is impossible to obtain such a bound using a weaker
$\psi_\alpha$ metric.

Unfortunately, if one is interested, as we are, in bounds on the
empirical process indexed by $F^2$ using complexity parameters of
$F$ itself, a contraction type argument only yields that
\begin{equation} \label{eq:emp-trivial-gamma-2-intro}
\E \sup_{f \in F} \left|\frac{1}{N} \sum_{i=1}^N f^2(X_i) - \E f^2
\right| \leq c (\sup_{f \in F} \|f\|_\infty)
\frac{\gamma_2(F,\psi_2)}{\sqrt{N}},
\end{equation}
which is unsatisfactory when dealing with a class of unbounded or
weakly bounded functions that only have nice tails. For such
classes, \eqref{eq:emp-trivial-gamma-2-intro} is meaningless.

An improvement to this contraction based estimate appeared in
\cite{KM} and later in \cite{MPT}, where is was shown that if $F$ is
a symmetric subset of the $L_2(\mu)$ unit sphere (i.e.
$\|f\|_{L_2(\mu)}=1$ and if $f \in F$ then $-f \in F$), one has
\begin{equation} \label{eq:emp-psi-2-gamma-2-intro}
\E \sup_{f \in F} \left|\frac{1}{N} \sum_{i=1}^N f^2(X_i) - \E f^2
\right| \leq c \max \left\{ (\sup_{f \in F} \|f\|_{\psi_2})
\frac{\gamma_2(F,\psi_2)}{\sqrt{N}}, \frac{\gamma_2^2(F,\psi_2)}{N}
\right\}.
\end{equation}
Thus, the diameter of $F$ in $L_\infty$ may be replaced by its
diameter in $\psi_2$.

There are many applications that follow from \eqref{eq:emp-psi-2-gamma-2-intro}.
For example, it was used in
\cite{MPT} to solve the approximate and exact reconstruction problems
(studied, e.g., in \cite{CRT,CT1,CT2}) in a rather general situation that includes any isotropic, subgaussian ensemble. However,
even \eqref{eq:emp-psi-2-gamma-2-intro} still leaves something to be
desired, since it too is meaningless for a large class of natural
measures. Indeed, consider the volume measure on an isotropic convex
body in $\R^n$, or more generally, an isotropic, log-concave measure
on $\R^n$. Again, if we set $F =\{\inr{x, \cdot} : x \in S^{n-1} \}$
-- the class of linear functionals of Euclidean norm one, it may
have a very bad diameter with respect to the $\psi_2(\mu)$ norm (as
bad as $\sqrt{n}$), whereas, thanks to Borell's inequality
(\cite{Bor}, see also \cite{MilSch}), its $\psi_1(\mu)$ diameter is
at most an absolute constant, independent of the dimension. Thus, it
seems natural to ask whether one may replace $d_{\psi_2}=\sup_{f \in
F}\|f\|_{\psi_2}$ in \eqref{eq:emp-psi-2-gamma-2-intro}, with
$d_{\psi_1}=\sup_{f \in F}\|f\|_{\psi_1}$. The main result of this
article is a positive answer to this question.

\noindent {\bf Theorem A.} {\it There exists an absolute constant
$c$ for which the following holds. If $F$ is a symmetric class of
mean-zero functions on $(\Omega,\mu)$ then
\begin{equation} \label{eq:emp-psi-1-gamma-2-intro}
\E \sup_{f \in F} \left|\frac{1}{N} \sum_{i=1}^N f^2(X_i) - \E f^2
\right| \leq c \max \left\{ d_{\psi_1}
\frac{\gamma_2(F,\psi_2)}{\sqrt{N}}, \frac{\gamma_2^2(F,\psi_2)}{N}
\right\},
\end{equation}
and a similar bound holds with high probability.}

\vskip 0.5cm

A key ingredient in the proof of Theorem A and our second main
result, deals with the structure of random coordinate projections of
a given class of functions that have nice tail properties. We will
be interested in the growth of the Euclidean norm of monotone
rearrangements of vectors in $P_\sigma F$: for every $1 \leq m \leq
N$ consider
$$
D_m=\sup_{v \in P_\sigma F} \left(\sum_{i=1}^m (v^2)_i^*
\right)^{1/2}=\sup_{f \in F} \max_{|I|=m} \left(\sum_{i \in I}
f^2(X_i) \right)^{1/2},
$$
where $(v_i^*)_{i=1}^N$ is a non-increasing rearrangement of
$(|v_i|)_{i=1}^N$. We will present high probability, sharp bounds
on the empirical diameters $D_m$ and use them in the proof of
Theorem A.

Let us consider a simple example that indicates which bound on $D_m$
one can hope for. Let $\mu$ be the canonical gaussian measure
on $\R^n$. Hence, if $(g_i)_{i=1}^n$ are independent, standard
normal random variables and $G=(g_1,...,g_n)$, then for every Borel
set $A \subset \R^n$, $\mu(A)=Pr(G \in A)$. Let $K \subset \R^n$,
consider $F=\{ \inr{x, \cdot} : x \in K\}$, the class of linear
functionals indexed by $K$, and put $(X_i)_{i=1}^N$ to be
independent, distributed according to $\mu$. Thus, $(X_i)_{i=1}^N$
are independent copies of $G$, and the coordinate projection of $F$
is given by $P_\sigma F = \left\{ ( \inr{x,X_i})_{i=1}^N \ : \ x \in
K\right\}$. Observe that there exists an absolute constant $c$
such that for every $1 \leq m \leq N$,
\begin{equation} \label{eq:lower-bound-intro}
\E \sup_{v \in P_\sigma F} \left(\sum_{i=1}^m (v^2)_i^*
\right)^{1/2} \geq c\left( \E \sup_{x \in K} \sum_{i=1}^n g_i x_i
+ \sup_{x \in K} \|x\|_{\ell_2^n} \sqrt{m\log(eN/m)} \right),
\end{equation}
where $\ell_2^n$ is the Euclidean norm on $\R^n$. Indeed, this lower
bound is evident because the first term is just the case $m=N=1$,
while the second term is an estimate for a single point $x \in K$
which has a maximal Euclidean norm.

The simple reasoning that leads to \eqref{eq:lower-bound-intro}
gives the impression that the estimate is far from sharp. However,
it turns out that there is an upper bound that holds in considerably
more general situations, and that matches the lower bound in the
gaussian case. The complexity parameter is, again, the $\gamma_2$
functional with respect to the $\psi_2$ norm, while the term that
represents the behavior of the ``worst" in the class is
$d_\alpha=\sup_{f \in F}\|f\|_{\psi_\alpha}$ for $1 \leq \alpha \leq
2$.

\vskip0.5cm

\noindent{\bf Theorem B.} For every $1 \leq \alpha \leq 2$ there is
a constant $c_\alpha$ that depends only on $\alpha$, and absolute
constants $c_1$ and $c_2$ for which the following holds. Let $F$ be
a class of mean-zero functions. Then, for every $u \geq c_1$, with
probability at least $1-\exp(-c_2 u \log N)$, for every $f \in F$
and every $1 \leq m \leq N$,
$$
\max_{|I|=m} \left(\sum_{i\in I} f^2(X_i) \right)^{1/2} \leq
c_\alpha u \left(\gamma_2(F,\psi_2) + d_{\psi_\alpha}m^{1/2}
\log^{1/\alpha}(eN/m) \right).
$$
\vskip0.5cm

To put this result in the right perspective let us return to the
gaussian example. If $\mu$ is the canonical gaussian measure on
$\R^n$ then the $\psi_2$ norm endowed on $\R^n$ is equivalent to the
Euclidean one. In particular, for every $m \leq N$, $\sup_{x \in
K}\|x\|_{\ell_2^n} \sqrt{m\log(eN/m)}$ and $\sup_{f \in
F}\|f\|_{\psi_2} \sqrt{m\log(eN/m)}$ are equivalent. Moreover, by
the Majorizing Measures Theorem (see \cite{Tal:book} and section
\ref{sec:pre}) and since the Euclidean and the $\psi_2$ metrics are
equivalent, so are $\E \sup_{x \in K} \sum_{i=1}^n g_i x_i$ and
$\gamma_2(F,\psi_2)$. Hence, the bound in Theorem B is sharp (up to
the absolute constants and the exact probabilistic estimate) for the
class of linear functionals indexed by a subset of $\R^n$ and with
respect to the gaussian measure.

Theorem B reveals useful information on the way vectors in $P_\sigma
F$ look like for a typical sample $(X_i)_{i=1}^N$. If $N$ is relatively
small, namely, when $d_{\psi_\alpha} \sqrt{N} \ll
\gamma_2(F,\psi_2)$, all the information one has is that the
Euclidean norm of any $P_\sigma f$ is at most of the order of
$\gamma_2(F,\psi_2)$. Then, for larger values of $N$ the situation
changes. For every $f \in F$ and
$$
\lambda = c_\alpha^\prime
d_{\psi_\alpha}
\log^{1/\alpha}(cNd_{\psi_\alpha}^2/\gamma_2^2(F,\psi_2)),
$$
the block $I=I(f)=\{i: |f(X_i)| \geq \lambda\}$ has small cardinality:
$$
\sup_{f \in F} |I(f)| \leq c_\alpha \gamma_2^2(F,\psi_2)/\lambda^2.
$$
Outside this
block, a monotone rearrangement of any $P_\sigma f$ is dominated
coordinate-wise by a rearrangement of the vector $(c_\alpha d_{\psi_\alpha}
\log^{1/\alpha} (eN/i))$.

Note that the behavior of the ``small" coordinates of each $P_\sigma
f$ is natural for a single $\psi_\alpha$ random variable. Indeed, it
is straightforward to verify that if $v$ is a $\psi_\alpha$ random
variable and $(v_i)_{i=1}^N$ is a vector of independent copies of
$v$, then with high probability, for every $i$, $v_i^* \leq
c\|v\|_{\psi_\alpha}\log^{1/\alpha} (eN/i)$. Thus, our results show
that for a random sample $\sigma$, the ``small coordinates" of any
$P_\sigma f$ are dominated by the typical behavior of
a sample of the function in the class with the maximal $\psi_\alpha$ norm.
From that point of view, each vector in $P_\sigma F$ can be decomposed into
a ``regular" part, which behaves as if $F$ has an envelope function
whose $\psi_\alpha$ norm is $d_{\psi_\alpha}$, and a ``peaky" part,
which is supported
on the block $I(f)$ and is bounded in $\ell_2^N$. The
blocks $I(f)$ take care of the possibility that vectors in $P_\sigma F$
have a few ``a-typical" large coordinates that are due to the
complexity of the whole class.

To formulate a weak version of the decomposition result (the full
one is presented in Theorem \ref{thm:decomposition}) we need two
preliminary definitions. First, for sets $A,B \subset \R^n$, $A+B =
\{a+b : a \in A, b \in B \}$ is the Minkowski sum of $A$ and $B$.
Second, we denote by $B_p^N$ the unit ball of $\ell_p^N=(\R^N,\| \
\|_p)$ and by $B_{\psi_\alpha^N}$ the unit ball of the $\psi_\alpha$
norm on $\R^N$, when viewed as the space of functions on the
probability space $\Omega=\{1,...,N\}$ endowed with the uniform
probability measure.

\noindent{\bf Theorem C.} {\it There exist absolute constants $c_1,...,c_6$
for which the following holds. Let $F$
be a class of mean-zero functions. For $1 \leq \alpha \leq 2$ and
for every $N$, set
$$
\lambda= c_1 \max\left\{  d_{\psi_\alpha}
\log^{1/\alpha}(c_2Nd_{\psi_\alpha}^2/\gamma_2^2(F,\psi_2)), 1
\right\}.
$$
Then, for every $t \geq c_3$, with probability at least
$1-2\exp(-c_4t\log N)$,
$$
P_\sigma F \subset c_5t \left(\gamma_2(F,\psi_2)B_2^N + (\lambda
B_\infty^N \cap c_6d_{\psi_\alpha^N} B_{\psi_\alpha^N}) \right).
$$
} 
\vskip 0.5cm
Theorem C extends and improves one of the main results from
\cite{Men:weak}. As we will explain in Section \ref{sec:deco},
it also extends an empirical processes version of a
theorem due to Rudelson on selector processes from \cite{Rud-selectors}.

Let us turn to the applications of the three theorems described
above that will be presented here. We will focus on properties of the
random operator $\Gamma = \sum_{i=1}^N \inr{X_i, \cdot} e_i$,
considered as an operator between an arbitrary $n$ dimensional
normed space $(\R^n, \| \ \|)$ and $\ell_p^N$, where $(X_i)_{i=1}^N$
are independent, distributed according to an isotropic, log-concave
measure $\mu$ on $\R^n$.

It is well known that many results in Asymptotic Geometric Analysis
have been obtained using certain specific random selection methods,
most often, according to the canonical gaussian measure on $\R^n$,
or with respect to the Haar measure on an appropriate Grassman
manifold. These selection methods, combined with analogs of Theorem
A and Theorem B for those models of randomness, lead to geometric
information on the structure of convex bodies, most notably, to
Dvoretzky type theorems and to low-$M^*$ estimates (see, e.g.
\cite{MilSch,Pis:book}).

We will show that sometimes it is possible to use more general
sampling methods and still obtain similar geometric results. In
particular, we will show that parts of the classical, gaussian based
theory, (e.g. ``standard shrinking" and low-$M^*$ estimates) may be
extended to log-concave ensembles. In fact, the gaussian parameter
associated with a convex body $K$, $\E \sup_{x \in K} \sum_{i=1}^n g_i x_i$,
which is used as a complexity parameter in the classical gaussian
based theory, is replaced in our results by $\gamma_2(K,\psi_2)$.
And, although the two complexity parameters are seemingly different,
it can be shown that they coincide if one resorts to the original
sampling methods.

Because of their general nature, Theorems A, B and C have many other
applications in very different directions, and these will not be
presented here. For example (out of many), our results can be used
to extend the analysis from the known cases to other ensembles
of the reconstruction problem, approximate
and exact (see, for example, \cite{CRT,CT1,CT2,MPT}), of the
statistical persistence problem \cite{GR,BMN} and of various
embedding problems. Some of
the applications are straightforward but others are more difficult,
since obtaining sharp estimates on the complexity parameter
$\gamma_2(F,\psi_2)$ can be nontrivial. To keep this article at a
reasonable length and to maintain its focus on the structural,
empirical processes oriented results, we chose to defer the
presentation of most of the applications to a later work.

The article is organized as follows. In Section \ref{sec:pre} we
will present preliminary results and several definitions we will
need. Then, in Section \ref{sec:diam} we will prove Theorem B and
Section \ref{sec:deco} will be devoted to the proof of Theorem C.
Theorem A will be proved in Section \ref{sec:conc}, and in Section
\ref{sec:applications} we will present some applications of the
three theorems.

\vskip1cm
\noindent{\bf Acknowledgments}

\noindent The research leading to these results has received funding
from the European Research Council under the European Community's
Seventh Framework Programme (FP7/2007-2013) / ERC grant agreement
n$^o$ [203134], from the Israel Science Foundation grant 666/06 and from the Australian Research Council grant DP0986563.
The author would like to thank M. Cwikel, M. Kozdoba and G. Lecu\'{e} for valuable discussions during the preparation of this article.

\section{Preliminaries} \label{sec:pre}
Let us begin with notational conventions. Throughout, all absolute
constants are positive numbers, denoted by $c,c_0,c_1,...$ etc.
Their value may change from line to line. We use
$\kappa_0,\kappa_1,...$ for constants whose value will remain
unchanged. By $A \sim B$ we mean that there are absolute constants
$c$ and $C$ such that $cB \leq A \leq CB$, and by $A \lesssim B$
that $A \leq CB$. For $1 \leq p \leq \infty$, $\ell_p^n$ is $\R^n$
endowed with the $\ell_p$ norm, which we denote by $\| \ \|_p$, and
$B_p^n$ is its unit ball. With a minor abuse of notation we denote
by $| \ |$ the cardinality of a set and the absolute value.

We say that $K \subset \R^n$ is a convex body if it is a compact,
convex and symmetric set (that is, if $x \in K$ then $-x \in K$)
with a nonempty interior. If $K$ is a convex body we denote by $\| \
\|_K$ the norm on $\R^n$ whose unit ball is $K$ and set $K^\circ=\{y
: \inr{x,y} \leq 1 \ \forall y \in K\}$ to be its polar body.

Given a probability measure $\mu$ and a sample $(X_i)_{i=1}^N$, we
will sometimes write $P_N f =N^{-1}\sum_{i=1}^N f(X_i)$ and $P f =
\E f$. Hence, the supremum of the empirical process indexed by $F$
is $\sup_{f \in F} |P_N f - Pf|$. Given $f \in F$ and $\sigma
\subset \{1,...,N\}$ we set $P_\sigma f = (f(X_i))_{i \in \sigma}$.

A significant part of our discussion will use basic properties of
sums of independent random variables that have nice tails. The
proofs of the claims presented here may be found, for example, in
\cite{LT}, \cite{GidlP} or \cite{VW}.

Recall that a random variable has a bounded $\psi_\alpha$ norm for
$1 \leq \alpha \leq 2$ if there is some constant $C$ for which
$\E \exp(|f|^\alpha/C^\alpha) \leq 2$, and in which case one sets
$$
\|f\|_{\psi_\alpha} = \inf\left\{ C : \E \exp(|f|^\alpha/C^\alpha) \leq 2 \right\}.
$$
One can show that  there is an absolute constant $c$ such
that if $f \in L_{\psi_\alpha}$, then for every $t \geq 1$,
$Pr \left(|f| \geq t\right) \leq
2\exp(-ct^\alpha/\|f\|_{\psi_\alpha}^\alpha)$.
Conversely, there is an absolute constant $c_1$ such that if $f$
displays a tail behavior dominated by $\exp(-t^\alpha/K^\alpha)$ for
some $1 \leq \alpha \leq 2$, then $f \in L_{\psi_\alpha}$ and
$\|f\|_{\psi_\alpha} \leq c_1K$.
We say that $X$ is a subgaussian random variable if $\|X\|_{\psi_2}
< \infty$.

\begin{Lemma} \label{lemma:sum-subgauss}
There exists an absolute constant $c$ for which the following holds.
Let $X$ be a mean-zero, subgaussian random variable and let
$X_1,...,X_k$ be independent copies of $X$. Then, for every fixed
$a=(a_1,...,a_k) \in \R^k$ $\|\sum_{i=1}^k a_i X_i \|_{\psi_2} \leq
c\|X\|_{\psi_2} \|a\|_2$. Thus, for every $t>0$,
$$
Pr \left( |\sum_{i=1}^k a_i X_i | \geq ct\|X\|_{\psi_2}\|a\|_2
\right) \leq 2\exp(-t^2/2).
$$
\end{Lemma}

In particular, if $(\eps_i)_{i=1}^N$ are independent, symmetric
$\{-1,1\}$-valued random variables, then for every
$(a_i)_{i=1}^N$,
$$
Pr \left( |\sum_{i=1}^N a_i \eps_i | \geq ct\|a\|_2 \right) \leq
2\exp(-t^2/2).
$$

For sums of independent $\psi_1$ random variables the situation is
more delicate, and one should expect two types of behaviors: an early
subgaussian decay followed by a subexponential one, as Bernstein's
inequality shows.

\begin{Lemma} \label{lemma:Bern}
There exists an absolute constant $c$ for which the following holds.
Let $X_1,...,X_N$ be independent copies of a mean-zero random
variable. Then, for any  $t>0$,
\begin{equation*}
\Pr \left(\left|\frac{1}{N}\sum_{i=i}^N X_i \right| >t \right)
\leq 2\exp \left(-c\, N\min \left( \frac{t}{\|X\|_{\psi_1}},\,
\frac{t^2}{\|X\|^2_{\psi_1}}\right)\right).
\end{equation*}
\end{Lemma}

This estimate may be extended to other values of $\alpha$. The next
lemma is a standard outcome of Corollaries 2.9 and 2.10 from
\cite{Tal:AJM94} (see \cite{ALPT2} for the proof).
\begin{Lemma} \label{lemma-sums-of-psi-alpha}
Let $1 \leq \alpha \leq 2$ and let $(X_i)_{i=1}^N$ be independent,
mean-zero random variables such that $\|X_i\|_{\psi_\alpha} \leq A$
for every $1 \leq i \leq N$. Then, for every $(a_i)_{i=1}^N \in
\R^N$ and any $t>0$,
$$
Pr\left(|\sum_{i=1}^N a_i X_i | \geq t A \right) \leq 2\exp\left(-c
\min\left\{\frac{t^2}{\|a\|^2_2},\frac{t^\alpha}{\|a\|^\alpha_{\alpha^*}}\right\}\right),
$$
where $1/\alpha+1/\alpha^* =1 $ and $c$ is an absolute constant.
\end{Lemma}

Next, let us turn to the main complexity parameter we will use -
Talagrand's $\gamma_2$ functional.
\begin{Definition} \label{def:gamma-2} \cite{Tal:book}
For a metric space $(T,d)$, an {\it admissible sequence} of $T$ is a
collection of subsets of $T$, $\{T_s : s \geq 0\}$, such that for
every $s \geq 1$, $|T_s| \leq 2^{2^s}$ and $|T_0|=1$. For $\beta \geq 1$,
define the $\gamma_\beta$ functional by
$$
\gamma_\beta(T,d) =\inf \sup_{t \in T} \sum_{s=0}^\infty
2^{s/\beta}d(t,T_s),
$$
where the infimum is taken with respect to all admissible sequences
of $T$. For an admissible sequence $(T_s)_{s \geq 0}$ we denote by
$\pi_s t$ a nearest point to $t$ in $T_s$ with respect to the metric
$d$.
\end{Definition}

When
considered for a set $T \subset L_2$, $\gamma_2$ has close connections with
properties of the canonical gaussian process indexed by $T$, and we
refer the reader to \cite{Dud-book,Tal:book} for detailed
expositions on these connections. One can show that under mild measurability
assumptions, if $\{G_t: t \in T\}$ is a centered gaussian process indexed by a set $T$ then
$$
c_1 \gamma_2(T,d) \leq \E \sup_{t \in T} G_t \leq c_2 \gamma_2(T,d),
$$
where $c_1$ and  $c_2$ are  absolute constants and for every $s,t
\in T$, $d^2 (s,t) = \E|G_s-G_t|^2$.  The upper bound is due to
Fernique \cite{F} and the lower bound is Talagrand's Majorizing
Measures Theorem \cite{Tal87}. Note that if $T \subset \R^n$,
$(g_i)_{i=1}^n$ are standard, independent gaussians and $G_t =
\sum_{i=1}^n g_i t_i$ then $d(s, t) = \|s-t\|_2$, and therefore
\begin{equation}
  \label{maj_meas}
  c_1 \gamma_2(T,\| \cdot \|_2) \leq  \E \sup_{t \in T} \sum_{i=1}^n
  g_it_i  \leq c_2 \gamma_2(T,\|\cdot \|_2).
\end{equation}

Note that a closely related complexity parameter that is used
to describe geometric properties of a convex body $K$ is
$$
M^*(K)=\int_{S^{n-1}} \|x\|_{K^\circ} d\sigma (x),
$$
where $\sigma$
is the Haar measure on the sphere $S^{n-1}$. This parameter is
gaussian in nature and it is straightforward to verify that
$$
\sqrt{n}
M^*(K) \sim \E \sup_{x \in K} \sum_{i=1}^k g_i x_i.
$$

It is well known that chaining methods lead to simple bounds on empirical processes.
Indeed, the following result is a combination of a chaining argument with Lemma
\ref{lemma:sum-subgauss} or with Lemma \ref{lemma:Bern}.

\begin{Theorem}  \label{thm:generic-chaining}
There exists an absolute constant $c$ for which the following holds.
If $F$ is a class of functions on $(\Omega,\mu)$ then for every
integer $N$,
\begin{equation*}
\E \sup_{f \in F} |P_N f-P f| \leq
c\left(\frac{\gamma_2(F,\psi_1)}{\sqrt{N}}+
\frac{\gamma_1(F,\psi_1)}{N} \right),
\end{equation*}
and
$$ \E \sup_{f \in F} |P_N f-P f| \leq
c\frac{\gamma_2(F,\psi_2)}{\sqrt{N}}.$$
Similar bounds hold with
high probability.
\end{Theorem}
Results of this flavor may be found in Chapters 1 and 2.7 of
\cite{Tal:book}.

Finally, in Section \ref{sec:applications} we will be interested
in isotropic, log-concave measures on $\R^n$.

\begin{Definition}
A symmetric probability measure $\mu$ on $\R^n$ is called isotropic
if for every $y \in \R^n$, $\int |\inr{x,y}|^2 d\mu(x) = \|y\|_2^2$.

We say that a measure $\mu$ on $\R^n$ is $L$-subgaussian if for
every $x \in \R^n$, $\|\inr{x,\cdot}\|_{\psi_2(\mu)} \leq
L\|\inr{x,\cdot}\|_{L_2(\mu)}$.

The measure $\mu$ is log-concave if for every $0<\lambda<1$ and
every nonempty Borel measurable sets $A,B \subset \R^n$,
$\mu(\lambda A+(1-\lambda)B) \geq \mu(A)^\lambda\mu(B)^{1-\lambda}$.
\end{Definition}
The canonical gaussian measure on $\R^n$ is clearly isotropic and
subgaussian, with $L$ being an absolute constant. Lemma
\ref{lemma:sum-subgauss} implies that the same holds for the uniform
measure on $\{-1,1\}^n$.

A typical example of a log-concave measure on $\R^n$ is the volume
measure of a convex body in $\R^n$, a fact that follows from the
Brunn-Minkowski inequality (see, e.g. \cite{Pis:book}). Moreover,
Borell's inequality \cite{Bor,MilSch} implies that there is an
absolute constant $c$ such that if $\mu$ is an isotropic,
log-concave measure on $\R^n$, then for every $x \in \R^n$,
$\|\inr{x,\cdot}\|_{\psi_1} \leq c\|\inr{x,\cdot}\|_{L_2} =
c\|x\|_2$. There are isotropic bodies with a subgaussian volume
measure -- for example, isotropic positions of $B_p^n$ for $p \geq
2$ \cite{BGMN}. However, the general situation is completely
different, and there are many examples of volume measures of
isotropic convex bodies in $\R^n$ for which linear functionals are
far from exhibiting a bounded $\psi_2$ behavior. In fact,
$\|\inr{x,\cdot}\|_{\psi_2}$ may be as large as $\sqrt{n} \|x\|_2$
(for example, $x=e_1$ and the volume measure on an
isotropic position of $B_1^n$). We refer the reader to
\cite{Gia:survey} for a survey on properties of the volume measure
of isotropic convex bodies and, more generally, of isotropic
log-concave measures on $\R^n$.

\section{Bounding the diameter} \label{sec:diam}

This section is devoted to the proof of Theorem B. Although we will present a
complete proof only for $\alpha=1$, we will indicate the very
minor modifications that are needed to prove it for any $1 \leq
\alpha \leq 2$.

The first step in the proof of Theorem B is to construct a good
cover of the Euclidean unit ball $B_2^N$, an idea which was used for
a very similar goal in the proof of the main result in \cite{ALPT}.
\begin{Definition}
If $A,B \subset \R^n$, we denote by $N(A,B)$ the smallest number of
points $x_i \in A$ such that $A \subset \bigcup_i (x_i +B)$.

If $\| \ \|$ is a norm on $\R^n$ and $B=\{x : \|x\| \leq \eps\}$,
then the set $\{x_i\}$ is called an $\eps$-cover of $A$ with respect
to the norm $\| \ \|$.
\end{Definition}
Clearly, if $B$ is an $\eps$ ball of some norm, then $N(A,B)$ is
the smallest cardinality of a set $\{x_i\}$ such that for every
$a \in A$, $\min_i \|a-x_i\| \leq \eps$.

Fix an integer $N$ and define the following sets: for $1 \leq \ell
\leq N/2$ put
$$
A_\ell = \left\{z \in B_2^N : \ |{\rm supp} (z)| \leq \ell, \
\|z\|_\infty \leq 1/\sqrt{\ell} \right\}.
$$

Let $\eps_\ell=\ell/N$, set $N_\ell \subset A_\ell$ to be an
$\eps_\ell$-cover of $A_\ell$ with respect to the $\ell_2^N$ norm
and let $P_I:\R^N \to \R^N$ be the orthogonal projection onto the
space spanned by the coordinates $(e_i)_{i \in I}$, that is,
$P_Ix=\sum_{i \in I} \inr{e_i, x }e_i$. A standard volumetric
estimate shows that for every convex body $K \subset \R^N$, $N(K,
\eps K) \leq (2/\eps)^n$. Therefore,
\begin{equation}
|N_\ell| \leq \sum_{|I|=\ell} N(P_I B_2^N,\eps_\ell B_2^N) \leq
\binom{N}{\ell} \left(\frac{2}{\eps_\ell}\right)^\ell \leq \exp(c_0
\ell \log(eN/\ell))
\end{equation}
for a suitable absolute constant $c_0$.

Fix an integer $m \leq N$ and assume that $m=2^{r_0}$ for some
integer $r_0$. Define the sets $B_m$ as follows:
\begin{equation} \label{eq:B_m}
B_m = \left\{ z \in B_2^N : \ |{\rm supp}(z)| \leq m, \ {\rm
supp}(z)=\bigcup_{r=0}^{r_0-1} I_r \ , \  P_{I_r} z \in N_{|I_r|}
\right\}
\end{equation}
where $(I_r)_{r=0}^{r_0-1}$ are disjoint sets of coordinates, $|I_0|=2$ and for $r
\geq 1$, $|I_r|=2^r$ (and thus the cardinality of their union is
$m$).

It is evident that $B_m$ consists of vectors in $B_2^N$ that can be
written as a sum over disjoint sets of coordinates $I_r$ of
cardinality $2^r$, and the projection onto each one of the ``blocks" $I_r$
belongs to the net $N_{|I_r|}$, and thus to $A_{|I_r|}$.

It is standard to verify that for every $m=2^{r_0}$,
$$
|B_m| \leq |N_2| \cdot \prod_{r=1}^{r_0-1} |N_{2^r}| \leq
\prod_{r=0}^{r_0-1} \exp(c_02^r \log (eN/2^r)) \leq \exp(c_1
m\log(eN/m)).
$$

Let $D_m = \sup_{f \in F} \sup_{|I|=m} \left(\sum_{i \in I} f^2(X_i)
\right)^{1/2}$. The next lemma shows that in order to bound $D_m$ it
is enough to consider the linearized process indexed by $F \times
B_m$ and defined by $(f,v) \to \sum_{i=1}^N f(X_i) v_i$. Although
Lemma \ref{lemma:approximation-by-B-sets} is a purely deterministic
result, it is formulated in the ``random" context in which it will
be used.

\begin{Lemma} \label{lemma:approximation-by-B-sets}
There exists an absolute constant $C$ such that for every $m \leq
N/2$ satisfying $m=2^{r_0}$ for some integer $r_0$, and for every
$X_1,...,X_N$,
$$
D_m \leq C \sup_{f \in F} \sup_{v \in B_m} \sum_{i=1}^N v_i f(X_i).
$$
\end{Lemma}
\proof Let $m=2^{r_0}$ for some integer $r_0$ and assume that $m
\leq N/2$. If $v \in B_2^N$ for which $|{\rm supp}(v)| \leq m$,
let $(v_i^*)_{i=1}^N$ be a monotone non-increasing rearrangement of
$(|v_i|)_{i=1}^N$ and put $v_{\sigma(j)}=v_j^*$, where $\sigma$ is
the suitable permutation of $\{1,...,N\}$. Consider the sets
$(I_r)_{r=0}^{r_0-1}$ defined as follows:
$I_0=\{\sigma(1),\sigma(2)\}$ are the largest two coordinate of
$(|v_i|)_{i=1}^N$, $I_1=\{\sigma(3),\sigma(4)\}$ are the two following that,
and so on -- $I_r = \{ \sigma(2^r+1),...,\sigma(2^{r+1})\}$ for $r
\geq 1$. Thus, $|I_r|=2^r$ for $r \geq 1$ and $|I_0|=2$. Since
$v_i^* \leq 1/\sqrt{i}$ then for every $j$,
$$
\|P_{I_j}v\|_\infty
\leq 1/(| \bigcup_{r < j} I_r|)^{1/2} =2^{-j/2}
$$
and thus $P_{I_j} v \in A_{2^j}$.
Let $\tilde{v} \in B_m$ be such that for every $1 \leq
r \leq r_0-1$, $P_{I_r} \tilde{v} \in N_{|I_r|}=N_{2^r}$,
$\|P_{I_r}v- P_{I_r} \tilde{v}\|_2 \leq 2^r/N$ and $\|P_{I_0}v-
P_{I_0} \tilde{v}\|_2 \leq 2/N$.

Therefore, $ \|v - \tilde{v}\|_2 \leq 2/N+ \sum_{r=1}^{r_0-1}
{2^{r}}/{N} \leq {m}/{N}$, and thus, if we set $U_m = \{v \in B_2^N
: |{\rm supp}(v)| \leq m\}$ then
\begin{align*}
D_m & =\sup_{f \in F, \ |I|=m} \left(\sum_{i \in I} f^2(X_i)
\right)^{1/2} = \sup_{f \in F, \ |I|=m} \sup_{v \in U_m} \sum_{i \in
I} v_i f(X_i)
\\
& \leq \sup_{f \in F, \ |I|=m} \sup_{v \in U_m} \sum_{i \in I} (v_i
-\tilde{v}_i)f(X_i) + \sup_{f \in F, \ |I|=m} \sup_{v \in U_m}
\sum_{i \in I} \tilde{v}_i f(X_i)
\\
& \leq (m/N) D_m + \sup_{f \in F} \sup_{v \in U_m} \sum_{i \in I}
\tilde{v}_i f(X_i)
\\
& \leq (m/N) D_m + \sup_{f \in F} \sup_{v \in B_m} \sum_{i=1}^N v_i
f(X_i).
\end{align*}
Since $m \leq N/2$ it is evident that $ D_m \leq 2 \sup_{f \in F}
\sup_{v \in B_m} \sum_{i=1}^N v_i f(X_i), $ as claimed.
\endproof

\begin{Remark}
Observe that for every $1 \leq j \leq r_0-1$ and every $v \in
B_{2^{j+1}}$, $P_{\bigcup_{r<j}I_r} v \in B_{2^{j}}$, a fact which
will be used in the dimension reduction procedure that is needed
in the proof of Theorem B.
\end{Remark}

We will need two simple observations about sums of centered random
variables, both of which follow from Lemma \ref{lemma:sum-subgauss}
and Lemma \ref{lemma-sums-of-psi-alpha}. First, if $\E f =0$ then
for every $t>0$ and any $I \subset \{1,...,N\}$,
$$
Pr \left( \left|\sum_{i \in I} v_i f(X_i) \right| \geq t \|P_I v\|_2
\|f\|_{\psi_2} \right) \leq 2\exp(-c_0t^2).
$$
Second, if $\E f =0 $ then for every $t>0$ and any $I \subset
\{1,...,N\}$,
\begin{align} \label{eq:Bernstein-type-eq}
Pr \left( \left|\frac{1}{|I|} \sum_{i \in I} v_i f(X_i) \right| \geq
t\|P_I v\|_{\infty} \|f\|_{\psi_1} \right) \leq
2\exp(-c_0|I|\min(t^2,t)),
\end{align}
where in both cases $c_0$ is an absolute constant.

Before proving Theorem B we need a few more definitions.
Let $E_\ell$ be the collection of
all subsets of $\{1,...,N\}$ of cardinality $\ell$. Note
that there is an absolute constant $\kappa_0$ such that for every
integer $1 \leq \ell \leq N$, $\exp(\kappa_0 \ell \log(eN/\ell))
\geq \max\{|E_\ell|,|B_\ell|\}$, and define $s_\ell$
to be the first integer which satisfies that $2^{2^{s}} \geq
\exp(\kappa_0 \ell \log(eN/\ell))$.

The chaining argument we will use for $ f \to \sup_{v \in B_m}
\sum_{i=1}^N v_i f(X_i) $ consists of three parts. First, when $s
\geq s_m$, the number of vectors in $B_m$ is much smaller than the
number of possible ``links" in all the chains, and thus no special
treatment is needed. In the middle part, when $s_2 \leq s < s_m$,
there will be a simultaneous reduction in the level $s$ and in the dimension
which will be achieved by passing from the set $B_m$ to the
sets $B_{m/2^{r}}$ for the correct value of $r$. Finally, when $s
\leq s_2$ no further chaining will be required because the
cardinality of the indexing sets is small enough.

Let us reformulate Theorem B.
\begin{Theorem} \label{thm:main-est}
For every $1 \leq \alpha \leq 2$ there are constant $c_\alpha$ and
$C_\alpha$ that depend only on $\alpha$, and there exist absolute
constants $c_1 \geq 1$ and $c_2$ for which the following holds. Let
$F$ be a class of mean-zero functions and let $(F_s)_{s \geq 0}$ be
an admissible sequence of $F$. Then, for every $t \geq c_1$ and
every integer $N$, with probability at least $1-2\exp(-c_2 t \log N
)$, for every $m \leq N$ and every $f \in F$,
$$
\sup_{v \in B_m} \sum_{i=1}^N v_i f(X_i) \leq c_\alpha t
\left(\sum_{s=0}^\infty 2^{s/2} \|\pi_s f -\pi_{s-1}f\|_{\psi_2} +
d_{\psi_\alpha}  \sqrt{m} \log^{1/\alpha}(eN/m) \right),
$$
where $d_{\psi_\alpha}=\sup_{f \in F}\|f\|_{\psi_\alpha}$.

In particular, with that probability, for every $m \leq N$,
$$
D_m \leq C_\alpha t \left(\gamma_2(F,\psi_2) + d_{\psi_\alpha}
\sqrt{m} \log^{1/\alpha}(eN/m) \right).
$$
\end{Theorem}
As we said, we will present the proof of Theorem \ref{thm:main-est}
only for $\alpha=1$. The proof for $1<\alpha \leq 2$ is identical,
with the exception that \eqref{eq:Bernstein-type-eq} is replaced by
an appropriate deviation estimate for $\psi_\alpha$ random
variables, as stated in Lemma \ref{lemma-sums-of-psi-alpha}.

\proof Let $\{F_s : s \geq 0\}$ be an admissible sequence of $F$ and
without loss of generality, assume that $m=2^{r_0}$ for some integer
$r_0$.

To begin the first part of the chaining argument, for every fixed
$f$ set $(\Delta_s f)_i=(\pi_{s} f-\pi_{s-1} f)(X_i)$. Then, for
every $f \in F$ and $v \in B_m$,
$$
\sum_{i=1}^N v_i f(X_i) = \sum_{s > s_m} \sum_{i=1}^N v_i (\Delta_s
f)_i + \sum_{i=1}^N v_i (\pi_{s_m}f)(X_i).
$$
Since the cardinality of the set $\Delta_s = \{\pi_{s} f-\pi_{s-1} f
: f \in F\}$ is at most $2^{2^{s+1}}$ and since $|B_m| \leq
\exp(\kappa_0 m \log (eN/m))$, then by the definition of $s_m$ and a
$\psi_2$ estimate, if $t$ is larger than an absolute constant, one
has
\begin{align*}
Pr & \left(\exists f \in F, \ v \in B_m : |\sum_{s > s_m}
\sum_{i=1}^N v_i (\Delta_s f)_i| \geq t\|v\|_2 \sum_{s > s_m}
2^{s/2}\|\Delta_s f\|_{\psi_2} \right)
\\
\leq &  |B_m| \cdot 2\sum_{s > s_m} |\Delta_{s}|\exp(-c_1 2^{s} t^2)
\leq 2\exp(-c_2 2^{s_m}t^2).
\end{align*}

Now, let us turn to the ``middle part", in which the structure of
vectors that belong to $B_m$ is used. First, consider the integers
$s_m$, $s_{m/2}$, etc. From the definition of $s_\ell$ it follows
that there is an absolute constant $c_3$ such that for every $1 \leq
\ell \leq N$, $s_\ell$ satisfies that
$$
2^{s_\ell} \geq c_3 \ell \log (eN/\ell), \ \ \ \ 2^{s_\ell-1} < c_3
\ell \log(eN/\ell).
$$
In particular, $2^{s_\ell-2} \leq c_3 (\ell/2) \log (eN/\ell) < c_3
(\ell/2) \log (eN/(\ell/2))$, implying that $s_\ell -1 \leq
s_{\ell/2} \leq s_\ell$. A similar argument shows that $s_{\ell/4} <
s_\ell$ if $\ell \leq N/2$, and thus, either $s_{\ell/2}=s_\ell-1$
or, if $s_{\ell/2}=s_\ell$, then $s_{\ell/4}=s_\ell-1$. In any case,
if one considers the sequence $s_\ell, s_{\ell/2}, ....,
s_{\ell/2^r}$, it decreases in steps of at most one and remains
constant on blocks of cardinality at most two.

Fix any $v \in B_m$, and one may assume that $|{\rm supp}(v)|=m$.
Let $\ell_r =m/2^r$ and put $I_{\ell_1} \subset {\rm supp}(v)$ to be a
set of $m/2$ coordinates such that $\|P_{I_{\ell_1}}v\|_\infty \leq
1/(m/2)^{1/2}$ (such a set of coordinates exists by the definition
of $B_m$). Denote by $J_1$ the complement of $I_{\ell_1}$ in ${\rm
supp}(v)$, and observe that $P_{J_1} v \in B_{m/2}$ (where, of course,
$I_{\ell_1}$ and $J_1$ depend on $v$). Hence,
\begin{equation*}
\sum_{i=1}^N v_i(\pi_{s_m}f)(X_i) = \sum_{i \in I_{\ell_1}}
v_i(\pi_{s_m}f)(X_i) + \sum_{i \in J_1} v_i(\pi_{s_m}f)(X_i),
\end{equation*}

\begin{equation} \label{eq:first-part}
\sum_{i \in I_{\ell_1}} v_i(\pi_{s_m}f)(X_i) = \sum_{i \in
I_{\ell_1}} v_i (\pi_{s_m}f-\pi_{s_{m/2}}f)(X_i) + \sum_{i \in
I_{\ell_1}} v_i (\pi_{s_{m/2}}f)(X_i),
\end{equation}
and
\begin{equation} \label{eq:second-part}
\sum_{i \in J_1} v_i(\pi_{s_m}f)(X_i) = \sum_{i \in J_1} v_i
(\pi_{s_m}f-\pi_{s_{m/2}}f)(X_i) + \sum_{i \in J_1} v_i
(\pi_{s_{m/2}}f)(X_i).
\end{equation}

We will estimate the first part of \eqref{eq:first-part} using a
$\psi_2$ argument and the second one using the $\psi_1$ information.
Indeed, there are at most $2^{2^{s_m+1}}$ elements of the form
$\pi_{s_m}f-\pi_{s_{m/2}}f$, and at most $|B_m|$ vectors $v$. Since
$|B_m| \leq \exp(\kappa_0m\log(eN/m))$ then from the definition of
$s_m$ it follows that with probability at least
$1-2\exp(-c_4t^22^{s_m})$, for every $f \in F$ and $v \in B_m$
\begin{equation} \label{eq:first-part-1}
|\sum_{i \in I_{\ell_1}} v_i (\pi_{s_m}f-\pi_{s_{m/2}}f)(X_i)| \leq
t2^{s_m/2} \|P_{I_{\ell_1}}v\|_2
\|\pi_{s_m}f-\pi_{s_{m/2}}f\|_{\psi_2}.
\end{equation}
To handle the second term, recall that every $v \in B_m$,
$\|P_{I_{\ell_1}}v\|_\infty \leq 1/(m/2)^{1/2}$ and
$|I_{\ell_1}|=m/2$. Hence, for every $f \in F$, $v \in B_m$ and
$u>0$,
$$
Pr \left( \left|\sum_{i \in I_{\ell_1}} v_i (\pi_{s_{m/2}}f)(X_i)
\right| \geq u\|\pi_{s_{m/2}}f\|_{\psi_1}
\frac{|I_{\ell_1}|}{(m/2)^{1/2}} \right) \leq
2\exp(-c_5|I_{\ell_1}|\min(u^2,u)).
$$
In particular, if one takes $u=t2^{s_m}/|I_{\ell_1}|$ (which is of
the order of $\log(eN/m)$), then by our estimates on the cardinality
of $B_m$ and the definition of $s_{m/2}$, it follows that with
probability at least $1-2\exp(-c_6t 2^{s_m})$, for every $f \in F$
and every $v \in B_m$,
$$
\left|\sum_{i \in I_{\ell_1}} v_i (\pi_{s_{m/2}}f)(X_i) \right| \leq t
d_{\psi_1} \sqrt{m} \log(eN/m).
$$

Turning to \eqref{eq:second-part}, the first term can be bounded
exactly as in \eqref{eq:first-part-1}, while in the second term of
\eqref{eq:second-part}, the required dimension reduction is
achieved:
all the vectors $P_{J_1} v$ belong to $B_{m/2}$ and the indexing
class is $F_{s_{m/2}}$.

The same argument can be repeated, by breaking each $J_1$ into
$I_{\ell_2}$ and its complement in $J_1$ (which we denote by
$J_2$), just as in \eqref{eq:first-part} and
\eqref{eq:second-part}. At the $r$-th step one begins with vectors
$P_{J_{r-1}}v$ that belong to $B_{m/2^{r-1}}$, and an indexing set
$F_{s_{m/2^{r-1}}}$. It follows that with probability at least $
1-4\exp(-c_6t^2 2^{s_{m/2^{r-1}}})-2\exp(-c_6t 2^{s_{m/2^{r-1}}}),
$ for every $f \in F$ and $v \in B_m$,
\begin{align*}
& \left|\sum_{i=1}^N (P_{J_{r-1}}v)_i
(\pi_{s_{m/2^{r-1}}}f)(X_i)\right|
\\
\leq & t 2^{s_{m/2^{r-1}}/2} \|\pi_{s_{m/2^{r-1}}}f
-\pi_{s_{m/2^{r}}}f\|_{\psi_2}
\left(\|P_{I_{\ell_r}}v\|_2+\|P_{J_r}v\|_2\right)
\\
+ & t d_{\psi_1}{2^{s_{m/2^{r-1}}}}\|P_{I_{\ell_r}}v\|_\infty
\\
+ & \left|\sum_{i=1}^N (P_{J_{r}}v)_i
(\pi_{s_{m/2^{r}}}f)(X_i)\right|.
\end{align*}
Since $2^{s_{m/2^{r}}} \sim (m/2^{r})\log(eN/(m/2^{r}))$ and
$\|P_{I_{\ell_r}} v\|_{\infty} \leq 1/(m/2^{r})^{1/2}$, then $
{2^{s_{m/2^{r-1}}}}\|P_{I_{\ell_r}}v\|_\infty \lesssim
(m/2^{r-1})^{1/2} \log(eN/(m/2^{r-1})). $ Moreover,
$\|P_{I_r}v\|_2+\|P_{J_r}v\|_2 \leq 2\|P_{J_{r-1}}v\|_2 \leq 2$,
and thus the first two terms are bounded by
\begin{align*}
& c_7t 2^{s_{m/2^{r-1}}/2} \|\pi_{s_{m/2^{r-1}}}f
-\pi_{s_{m/2^{r}}}f\|_{\psi_2}
\\
+ & c_7 t d_{\psi_1} (m/2^{r-1})^{1/2}
\log\left(eN/(m/2^{r-1})\right).
\end{align*}
Hence, if we continue in this fashion until $s_2=s_{m/2^{r_0-1}}$,
it follows that for $t \geq c_8$, with probability at least
\begin{equation} \label{eq:probab-estimate}
1-4\sum_{r=1}^{r_0} \left( \exp(-c_6t^2 2^{s_{m/2^{r-1}}}) +
\exp(-c_6t2^{s_{m/2^{r-1}}})\right),
\end{equation}
for every $f \in F$ and every $v \in B_m$,
\begin{align*}
& \left|\sum_{i=1}^N v_i (\pi_{s_m}f)(X_i)\right|
\\
\leq & c_{9}t\left( \sum_{r=1}^{r_0} \left(2^{s_{m/2^r}/2}
\|\pi_{s_{m/2^{r-1}}}f -\pi_{s_{m/2^{r}}}f\|_{\psi_2}\right) +
d_{\psi_1} \sqrt{m}\log(eN/m) \right)
\\
+ & \left|\sum_{i=1}^N (P_{J_{r_0-1}}v)_i (\pi_{s_2}f)(X_i) \right|.
\end{align*}
Observe that the elements of the sequence $(s_{m/2^{r}})_{r=1}^{r_0}$ 
belong to the
interval $[s_1,s_m]$. Also, this sequence decreases in steps of at
most one, and each integer is repeated at most twice. Hence,
$$
\sum_{r=1}^{r_0} 2^{s_{m/2^r}/2} \|\pi_{s_{m/2^{r-1}}}f
-\pi_{s_{m/2^{r}}}f\|_{\psi_2} \leq 2 \sum_{s=s_2}^{s_m}
2^{s/2}\|\Delta_s(f)\|_{\psi_2},
$$
and the probabilistic estimate in \eqref{eq:probab-estimate} is at
least $ 1-2\exp(-c_{10}2^{s_2} t) \geq 1-2\exp(-c_{11} t \log{N}),
$ because $t \geq c_8$.

Finally, for the last step, consider the sets supported on at most two
coordinate, and thus $\log |F_{s_2}|, \ \log |B_2| \lesssim \log N$.
Therefore, with probability at least $1-2\exp(-c_{12} t\log{N})$,
for every $f \in F$ and $v \in B_m$,
$$
\left|\sum_{i=1}^N (P_{J_{r_0-1}}v)_i (\pi_{s_2}f)(X_i)\right| \leq
c_{13}td_{\psi_1} \log{N}.
$$
Summing the three parts, if follows that for every $t \geq C_0$ and for
every $m \leq N$, with probability at least $1-C_1\exp(-C_2 t
\log{N})$, for every $f \in F$ and every $v \in B_m$,
$$
\left|\sum_{i=1}^N v_i f(X_i) \right| \leq C_3t
\left(\sum_{s=1}^\infty 2^{s/2} \|\Delta_s(f)\|_{\psi_2} +
d_{\psi_1} \sqrt{m} \log(eN/m) \right).
$$
Since there are at most $N$ possible values of $m$, the same
holds for all $m \leq N$ uniformly, as claimed.
\endproof

Theorem B can be extended to other $\ell_p$ norms. Indeed, for $1
\leq p <2$ and any $I \subset \{1,...,N\}$, $\|x\|_{\ell_p^I} \leq
|I|^{1/p-1/2} \|x\|_2$. Hence, with probability at least
$1-2\exp(-c_1 t \log{N})$, for every $f \in F$ and $I \subset
\{1,...,N\}$,
\begin{equation} \label{eq:diameter-ell_p-small}
\left(\sum_{i \in I} |f|^p(X_i) \right)^{1/p} \leq c_\alpha t
\left(\gamma_2(F,\psi_2) |I|^{1/p-1/2} + d_{\psi_\alpha}|I|^{1/p}
\log^{1/\alpha}(eN/|I|) \right).
\end{equation}
For $p>2$ let $m_0$ be the smallest integer for which
$$
\gamma_2(F,\psi_2) \leq d_{\psi_\alpha} \sqrt{m}
\log^{1/\alpha}(eN/m).
$$
Then, by Theorem B, for every $|I| < m_0$,
$$
\left(\sum_{i \in I} |f|^p(X_i) \right)^{1/p} \leq \left(\sum_{i \in
I} |f|^2(X_i) \right)^{1/2} \leq 2c_\alpha t \gamma_2(F,\psi_2).
$$
For larger values of $|I|$, if we denote $(u_i)_{i=1}^N =
(f(X_i))_{i=1}^N$ then for $j \geq m_0$,
$$
u_j^* \leq \left(\frac{1}{j} \sum_{i=1}^j (u^2)_i^* \right)^{1/2}
\leq 2c_\alpha t d_{\psi_\alpha} \log^{1/\alpha} (eN/j).
$$
Hence, by the triangle inequality, for $|I| \geq m_0$,
\begin{equation} \label{eq:diameter-ell_p-large}
\left(\sum_{i \in I} |f|^p(X_i) \right)^{1/p} \leq c_{\alpha,p} t
\left(\gamma_2(F,\psi_2) + d_{\psi_\alpha}|I|^{1/p} \log^{1/\alpha}
(eN/|I|) \right).
\end{equation}
Let us mention that the estimate for $p=1$ (the weakest of
all the estimates for $1 \leq p \leq 2$) was proved in
\cite{Men:weak} using a simpler chaining argument.

\subsection{Optimality} \label{sec:diam-optimality}
We begin this section by recalling the observation made in the
introduction, that Theorem B is sharp when $F$ is a class of linear
functionals on $\R^n$ and $\mu$ is the canonical gaussian measure on
$\R^n$:
\begin{Lemma} \label{lower-est-gaussian}
There exists an absolute constant $c$ for which the following holds.
Let $K \subset \R^n$, set $F=\{\inr{x,\cdot} : x \in K\}$ and put
$\mu$ to be the canonical gaussian measure on $\R^n$. Then, for
every integer $N$ and any $1 \leq m \leq N$,
$$
\E \sup_{f \in F} \max_{|I|=m} \left( \sum_{i \in I} f^2(X_i)
\right)^{1/2} \geq c \left( \gamma_2(F,\psi_2) + d_{\psi_2}
\sqrt{m\log(eN/m)}\right).
$$
\end{Lemma}

Although Lemma \ref{lower-est-gaussian} indicates that Theorem B cannot
be improved, one might argue that it is a somewhat degenerate
case, because of the equivalence between the $\psi_2$ norm and the
$L_2$ one. The next lemma shows that in general, one cannot
replace the $\psi_2$ norm in the $\gamma_2$ term with any other
$\psi_\alpha$ norm for $\alpha <2$.
\begin{Lemma} \label{lemma:psi-alpha-lower-bound}
There exists an absolute constant $c_1$ for which the following
holds. For every integer $N$, $1 \leq \alpha <2$ and a number $R$,
there is a probability space $(\Omega,\mu)$ and a class $F$
consisting of mean-zero functions on $(\Omega,\mu)$, such that if
$(X_i)_{i=1}^N$ are independent, distributed according to $\mu$,
then with probability at least $c_1$,
$$
\sup_{f \in F} \left|\sum_{i=1}^N f(X_i) \right| \geq R
\gamma_2(F,\psi_\alpha)\sqrt{N},
$$
and in particular, $\sup_{f \in F} \left(\sum_{i=1}^N f^2(X_i)
\right)^{1/2} \geq R \gamma_2(F,\psi_\alpha)$.

\end{Lemma}

\begin{Remark} \label{rem:psi-alpha-lower-bound}
As we indicated in the introduction, Lemma
\ref{lemma:psi-alpha-lower-bound} shows that in general, $\E \sup_{f
\in F} |N^{-1}\sum_{i=1}^N f(X_i)-\E f |$ cannot be controlled using
a weaker deterministic parameter than $\gamma_2(F,\psi_2)/\sqrt{N}$.
\end{Remark}

For the proof of Lemma \ref{lemma:psi-alpha-lower-bound} we need the
following formulation of the Paley-Zygmund inequality \cite{GidlP}.
\begin{Lemma} \label{lemma:Pa-Zy}
Let $Z$ be a random variable. Then, for every $q > p \geq 1$ and $0
< \lambda <1$,
$$
Pr\left(|Z| \geq \lambda \|Z\|_{L_p} \right) \geq \left(
\left(1-\lambda^p \right) \left(\|Z\|_{L_p}/\|Z\|_{L_q}\right)^p
\right)^{q/(q-p)}.
$$
\end{Lemma}

\noindent {\bf Proof of Lemma \ref{lemma:psi-alpha-lower-bound}.}
Fix $1 \leq \alpha <2$ and an integer $n$. Let $Y$ be a symmetric
random variable with density $c_\alpha \exp(-|t|^\alpha)$ and set
$X=(Y_1,...,Y_n) \in \R^n$, a vector of independent copies of $Y$.
Consider the probability space $(\R^n,\mu)$, with $\mu$ defined by
$\mu(A)=Pr(X \in A)$, let $(e_i)_{i=1}^n$ be the standard basis of
$\R^n$, set $K=\{e_i/\sqrt{\log(i+1)} : 1 \leq i \leq n\}$ and put
$$
F=\{\inr{e_i,\cdot}/\sqrt{\log(i+1)} : 1 \leq i \leq n\}.
$$

One can show (see, for example, Proposition 7 in \cite{BGMN}) that
if $(x_i)_{i=1}^n$ is nonnegative and non-increasing, then for every
$p \geq 1$,
\begin{equation} \label{eq:optimality-moments}
\|\sum_{i=1}^n x_i Y_i\|_{L_p} \sim p^{1/\alpha} \|(x_i)_{ i \leq p}
\|_{\alpha^*} + \sqrt{p} \|(x_i)_{i >p}\|_2,
\end{equation}
where $\| \ \|_{\alpha^*}$ is the $\ell_{\alpha^*}$ norm for
$\alpha^*$ satisfying $1/\alpha+1/\alpha^*=1$. Since $\alpha \leq 2$
then $\alpha^* \geq 2$ and thus $\|\sum_{i=1}^n x_i Y_i\|_{L_p} \leq
c_1p^{1/\alpha}\|x\|_2$. In particular, for every $x \in \R^n$, $
\|\inr{x, \cdot}\|_{\psi_\alpha} \leq c_2 \|x\|_2. $ Moreover,
$\|\inr{x, \cdot}\|_{\psi_\alpha} \geq
c_3\|x\|_2$, implying that the $\ell_2^n$ and the $\psi_\alpha(\mu)$
norms are equivalent on $\R^n$.

It is also straightforward to show
that there is an absolute constant $c_4$ such that if
$(g_i)_{i=1}^\infty$ are independent, standard gaussian variables
then for every $m$,
$$
\E \max_{1 \leq i \leq m} \frac{g_i}{\sqrt{\log(i+1)}} \leq c_4.
$$
Therefore, by the Majorizing Measures Theorem
$$
\gamma_2(F,\psi_\alpha) \leq c_2 \gamma_2(K, \ell_2^n)  \leq  c_5
\E \max_{1 \leq i \leq n} \frac{g_i}{\sqrt{\log{(i+1)}}} \leq c_6.
$$

On the other hand, fix $N$ to be named later and consider $q > p
\geq N$. Observe that by \eqref{eq:optimality-moments}, for these
values of $q,p$ and $N$, if $(Y_i)_{i=1}^N$ are independent copies
of $Y$ then
$$
\left\|\sum_{i=1}^N Y_i \right\|_{L_p} \sim
p^{1/\alpha} N^{1-1/\alpha} \ \  {\rm and} \ \ \left\|\sum_{i=1}^N Y_i
\right\|_{L_q} \sim q^{1/\alpha} N^{1-1/\alpha}.
$$
Let $X_1,...,X_N$ be independent copies of the random vector $X$,
set $Y_{i,j}$ to be the $j$-th coordinate of $X_i$ and put
$Z_j=\sum_{i=1}^N Y_{i,j}$. Applying the Paley-Zygmund inequality,
it follows that there are $\beta>1$ and $c_7$, both depend on $\alpha$, such
that if $p=c_7 \log n$ and $q =\beta p$, then for every $j$,
$$
Pr
\left(|Z_j| \geq c_8(\log^{1/\alpha}n) N^{1-1/\alpha}/2 \right) =
Pr \left(|Z_j| \geq \|Z_j\|_{L_p}/2 \right) \geq
1/n.
$$
Hence, by the independence of $(Z_j)_{j=1}^n$,
$$
Pr \left( \exists 1 \leq j \leq n, \ |Z_j|  \geq c_8 N^{1-1/\alpha}
\log^{1/\alpha}n \right) \geq c_9.
$$
In particular, with that probability,
\begin{align*}
\sup_{f \in F } \left|\sum_{i=1}^N f(X_i) \right| =& \max_{1 \leq j
\leq n} \left| \sum_{i=1}^N \inr{\frac{e_j}{\sqrt{\log{(j+1)}}},X_i}
\right| = \max_{1 \leq j \leq n} \frac{1}{\sqrt{\log(j+1)}} \left|
\sum_{i=1}^N Y_{i,j} \right|
\\
\geq & c_{10} \frac{N^{1-1/\alpha} \log^{1/\alpha}n
}{\sqrt{\log(n+1)}} \gtrsim  \sqrt{N} \left(\frac{\log n}{N}
\right)^{1/\alpha -1/2}.
\end{align*}
All that remains now is to find the connection between $N$ and $n$,
where we already assumed that $p =c_7 \log n \geq N$. Clearly,
if $N \ll \log n$ then $\left(\log n /N\right)^{1/\alpha -1/2}$ can
be made to be arbitrarily large by increasing $n$, as claimed.
\endproof

\section{Decomposing $F$} \label{sec:deco}
Here, we will present a decomposition of $F$ into the sum of two
sets, representing its peaky and regular parts. We will
show that for every $N$, one can truncate functions in $F$ at the
level
$$
\lambda \sim  d_{\psi_\alpha} \log^{1/\alpha}(cd_{\psi_\alpha}
N^{1/2}/\gamma_2(F,\psi_2)).
$$
The resulting unbounded or peaky part of each $f
\in F$ has coordinate projections with a well behaved $\ell_2^N$
norm and short support. On the other hand, the regular part of $f$
is bounded in $L_\infty$ by $\lambda$, and, moreover, its typical coordinate
projection is contained in $cd_{\psi_\alpha} B_{\psi_\alpha^N}$.
Thus, the regular part of $F$ behaves as if $F$ has an envelope
function $W(x)=\sup_{f \in F}|f(x)|$ with $\|W\|_{\psi_\alpha} \leq d_\alpha$.

This decomposition gives a hint of why it is reasonable to hope that
the supremum of the empirical process $\sup_{f \in F} \left|P_N f^2 - P f^2 \right|$
is well behaved. Although the peaky part of $F$ exhibits no
concentration, its $\ell_2^N$ diameter is small, and thus there is
no need for cancelation to control it. Since the regular part of
$F$ behaves as if it has a reasonable envelope function, powers
concentrate around their mean uniformly.

To formulate the decomposition theorem (which implies Theorem C) we
will use the following observations. Recall that if $x \in \R^N$
then for $1 \leq \alpha \leq 2$, $\|x\|_{\psi_\alpha^N}=\inf\{C:
N^{-1}\sum_{i=1}^N \exp((|x_i|/C)^\alpha) \leq 2\}$. It follows that
for every $x \in \R^N$,
$$
x_i^* \leq c \|x\|_{\psi_\alpha^N}
\log^{1/\alpha} (eN/i),
$$
and, in fact, this behavior of a monotone
rearrangement characterizes the $\psi_\alpha^N$ norm. It is also
standard to verify that if $X$ is $\psi_\alpha$ random variable on
$(\Omega,\mu)$ and $(X_i)_{i=1}^N$ are independent copies of $X$,
then for $t \geq c_0$, with probability at least $1-2\exp(-t^\alpha \log N)$, for
every $i$,
$$
X_i^* \leq c_1t\|X\|_{\psi_\alpha} \log^{1/\alpha}(eN/i).
$$
Hence, a typical coordinate projection of an independent sample of a single
function $f \in L_{\psi_\alpha}$ satisfies that with high probability,
$\|(f(X_i))_{i=1}^N\|_{\psi_\alpha^N} \leq c_2 \|f\|_{\psi_\alpha}$.
In what follows, given $v=(f(X_i))_{i=1}^N$ we will sometimes denote the random
norms $ \|f\|_{\ell_p^N}$ and $\|f\|_{\psi_\alpha^N}$ by
$\|v\|_{p}$ and $\|v\|_{\psi_\alpha^N}$ respectively.

\begin{Theorem} \label{thm:decomposition}
There exist absolute constants $c_0,...,c_7$ for which the
following holds. For any $1 \leq \alpha \leq 2$ and an integer $N$
set
$$
\lambda = c_0 d_{\psi_\alpha} \max\left\{
\log^{1/\alpha}\left(c_0d_{\psi_\alpha}^2
N/\gamma_2^2(F,\psi_2)\right), 1\right\}.
$$
For any $t \geq c_1$ there are sets $F_1$ and $F_2$ that depend on
$N$, $\lambda$ and $t$ such that $F \subset F_1 + F_2$, and with
probability at least $1-2\exp(-c_2t\log{N})$,
\begin{description}
\item{1.}  $\sup_{f \in F_1} \|f\|_{\ell_2^N} \leq c_3 t \gamma_2(F,\psi_2)$, $\sup_{f \in
F_1} |{\rm supp} (P_\sigma f) | \leq c_3
\gamma_2^2(F,\psi_2)/\lambda^2$, and $\sup_{f \in F_1} \E |f|^2 \leq
c_3 \gamma_2^2(F,\psi_2)/N$.
\item{2.} $\sup_{f \in F_2} \|f\|_{L_\infty} \leq \lambda t$ and
$\sup_{f \in F_2} \|f\|_{\psi_\alpha^N} \leq c_4 t d_{\psi_\alpha}$.
\item{3.} For every $u \geq c_5$, with probability at
least $1-2\exp(-c_6u^2)$, one has $ \sup_{f \in F_2} \left|P_N f^2 -
P f^2 \right| \leq c_7u t\lambda {\gamma_2(F,\psi_2)}/{\sqrt{N}}$ .
\end{description}
\end{Theorem}

The proof of Theorem \ref{thm:decomposition} requires all the
information we have about the structure of the set
$P_\sigma F = \left\{ \left(f(X_i)\right)_{i=1}^N : f \in F \right\}
\subset \R^N$. Our starting point is the next observation.
\begin{Lemma} \label{lemma:coordinates}
There exist absolute constants $c_1$, $c_2$ and $c_3$ for which the
following holds. Let $v \in \R^N$ for which there are $A$, $B$ and
$1 \leq \alpha \leq 2$ such that for every $I \subset \{1,...,N\}$,
$$
\left(\sum_{i \in I} v_i^2 \right)^{1/2} \leq A + B
\sqrt{|I|}\log^{1/\alpha}\left(eN/|I|\right).
$$
If $\beta \geq c_1B\max\{ \log^{1/\alpha}(c_2NB^2/A^2),1\}$ and
$E_\beta = \{i: |v_i| \geq \beta \}$, then
$$
|E_\beta| \leq \max \left\{ \frac{4A^2}{\beta^2}, eN
\exp(-(\beta/2B)^{\alpha})\right\} \ \ {\rm and} \ \ (\sum_{i \in
E_\beta} v_i^2)^{1/2} \leq c_3 A.
$$
\end{Lemma}

\proof Clearly, for every integer $n$, $\|x\|_{\ell_1^n} \leq
\sqrt{n} \|x\|_{\ell_2^n}$. Hence, for every $I \subset
\{1,...,N\}$, $ \sum_{i \in I} |v_i| \leq A \sqrt{|I|} + B
|I|\log^{1/\alpha}(eN/|I|). $ Let $E_\beta = \{i: |v_i| \geq \beta
\}$ and note that
\begin{equation*}
\beta |E_\beta| \leq \sum_{i \in E_\beta} |v_i| \leq A
|E_\beta|^{1/2} + B |E_\beta|\log^{1/\alpha}(eN/|E_\beta|).
\end{equation*}
If $B |E_\beta|\log^{1/\alpha}(eN/|E_\beta|) \leq \beta |E_\beta|/2$
then $|E_\beta| \leq 4A^2/\beta^2$. Otherwise, if the reverse
inequality holds, then $|E_\beta| \leq eN
\exp(-(\beta/2B)^{\alpha})$. Thus,
\begin{equation} \label{eq:tail-estimate}
|\{i: |v_i| \geq \beta\}| \leq \max \left\{ \frac{4A^2}{\beta^2}, eN
\exp(-(\beta/2B)^{\alpha})\right\}.
\end{equation}
To complete the proof, let $ \beta \geq c_1B
\max\{\log^{1/\alpha}(c_2NB^2/A^2),1\}$. Therefore, $|E_\beta| \leq
{A^2}/{\beta^2}$, and thus, for our choice of $\beta$,
$$
\left(\sum_{i \in E_\beta} v_i^2 \right)^{1/2} \leq A + B
\frac{A}{\beta} \log^{1/\alpha}\left(eN\beta^2/A^2\right) \leq
c_3A.
$$
\endproof

\noindent {\bf Proof of Theorem \ref{thm:decomposition}.} Fix $t
\geq c_0$ and recall that by Theorem B, with probability at least
$1-2\exp(-c_1t\log n)$, for every $v \in P_\sigma F$ the assumptions
of Lemma \ref{lemma:coordinates} hold with $A \sim t
\gamma_2(F,\psi_2)$ and $B \sim t d_{\psi_\alpha}$. Just as in Lemma
\ref{lemma:coordinates}, set
$$
\beta \sim B\max\{
\log^{1/\alpha}(c_2NB^2/A^2),1\} \equiv \lambda t.
$$
Let $\phi(f)={\rm sgn}(f) \min \{|f|,\beta\}$ and
$\psi(f)=f-\phi(f)$, put $ F_1=\{\psi(f) : f \in F\}$, $F_2 =
\{\phi(f) : f \in F\}$ and observe that $F \subset F_1+F_2$.

Let us consider $F_1$, which is the unbounded part of $F$. Note that
for every $f \in F$, if we set $u_i=\left(\psi(f)\right)(X_i)$, then
$ \{ i: |u_i| \geq \beta \} \subset \{i: |f(X_i)| \geq \beta \} =
E_\beta, $ and on that set $|u_i| = |f(X_i)|-\beta$. Hence, by Lemma
\ref{lemma:coordinates},
$$
\left(\sum_{i \in E_\beta}
u_i^2\right)^{1/2} \leq c_3t \gamma_2(F,\psi_2).
$$
Also, since
$\|f\|_{\psi_\alpha} \leq d_{\psi_\alpha}$ then by integrating the
tail, one may verify that
$$
\E |\psi(f)|^2 \leq c_4 \beta^2
\exp(-c_5 (\beta/d_{\psi_\alpha})^{\alpha}) \leq
c_6\gamma_2^2(F,\psi_2)/N,
$$
proving the first part of the claim.

Turning to the second part, note that if $f \in F$ and $w=P_\sigma
 (\phi(f))$ then $\|w\|_\infty \leq \beta$. Let $m=A^2/\beta^2$ and observe that
$$
\gamma_2(F,\psi_2) \sim
d_{\psi_\alpha}  \sqrt{m} \max\{ \log^{1/\alpha}(eN/m),1\}.
$$
First,
assume that $m \leq N$. Therefore, since $\beta \leq c_7t
d_{\psi_\alpha} \log^{1/\alpha}(eN/m)$ then for every $j \leq m$, $
w_j^* \leq \beta \leq c_7t d_{\psi_\alpha} \log^{1/\alpha}(eN/j)$.
Moreover, if $j \geq m$ then by Theorem B,
\begin{align*}
\left(\sum_{i \leq j} (w^2)_i^*\right)^{1/2} & \lesssim t\left(\gamma_2(F,\psi_2)
+ d_{\psi_\alpha}\sqrt{j} \log^{1/\alpha}(eN/j)\right)
\\
& \lesssim td_{\psi_\alpha}\sqrt{j} \log^{1/\alpha}(eN/j).
\end{align*}
Therefore,
$$
w_j^* \leq
\left(\frac{1}{j}\sum_{i=1}^j (w^2)_i^* \right)^{1/2} \leq c_8t
d_{\psi_\alpha} \log^{1/\alpha}(eN/j),
$$
and thus, $\sup_{f \in
F_2}\|f\|_{\psi_\alpha^N} \lesssim t d_{\psi_\alpha}$, as claimed.

On the other hand, if $m \geq N$ then $\lambda \sim
d_{\psi_\alpha}$. Hence, $\sup_{f \in F_2} \|f\|_{L_\infty} \leq
\beta \lesssim td_{\psi_\alpha}$, implying that $\sup_{f \in
F_2}\|f\|_{\psi_\alpha^N} \lesssim td_{\psi_\alpha}$.

It remains to estimate the supremum of the empirical process indexed
by $|\phi(f)|^2$. Since $\phi(x)={\rm sgn}(x) \min\{x,\beta\}$ is
$1$-Lipschitz, then for every $f_1 , f_2 \in F$, $\left|
\left|\phi(f_1)\right|^2 - \left|\phi(f_2)\right|^2 \right| \leq
2\beta |f_1-f_2|$ pointwise. In particular, $ \left\|
\left|\phi(f_1)\right|^2 - \left|\phi(f_2)\right|^2
\right\|_{\psi_2} \leq 2\beta \|f_1 - f_2\|_{\psi_2}. $ Therefore,
by a standard chaining argument, for every $u \geq c_9$, with
probability at least $1-2\exp(-c_{10}u^2)$, $ \sup_{h \in F_2}
\left|P_N h^2 - P h^2 \right| \leq c_{11}u
\beta{\gamma_2(F,\psi_2)}/{\sqrt{N}}$, as claimed.
\endproof

Theorem \ref{thm:decomposition} should be compared with a result due to Rudelson
(see \cite{Rud-selectors}). Although Rudelson's result was formulated for selector
processes, an analogous result holds for empirical processes and with essentially
the same proof as the original one.

\begin{Theorem} \label{thm:Rud}
For every $0<\delta<1$ there is a constant $c(\delta)$ for which the following holds. If $F$ is a class of mean-zero functions then there are sets $F_1$ and $F_2$ such that $F \subset F_1+F_2$, and with $\mu^N$-probability at least $1-\delta$,
$$
\sup_{f \in F_1} \|f\|_{\ell_2^N} \leq c(\delta)\sqrt{N}d_{L_2}, \ \ \ \sup_{f \in F_2} \|f\|_{\ell_1^N} \leq c(\delta) \E \sup_{f \in F} |\sum_{i=1}^N \eps_i f(X_i) |.
$$
In particular, with probability at least $1-\delta$,
$$
P_\sigma F \subset c(\delta) \left( R_N \sqrt{N}B_1^N + d_{L_2} \sqrt{N} B_2^N\right),
$$
where
$$
R_N = \frac{1}{\sqrt{N}} \E \sup_{f \in F} |\sum_{i=1}^N \eps_i f(X_i) | \ \ {\rm and} \ \  d_{L_2}=\sup_{f \in F} (\E f^2)^{1/2}.
$$
\end{Theorem}

Let us compare Theorem \ref{thm:decomposition} with Theorem \ref{thm:Rud}. First, observe that the first two parts of Theorem \ref{thm:decomposition} imply that
$$
P_\sigma F \subset c(\delta)\left(\gamma_2(F,\psi_2)B_2^N+\lambda B_\infty^N \cap cd_{\psi_\alpha}B_{\psi_\alpha^N}\right),
$$
and that for large values of $N$, that is, when
$\lambda \geq 1$, it is evident that $\gamma_2(F,\psi_2) \leq \sqrt{N}d_{\psi_\alpha}$. In particular, since $B_{\psi_\alpha^N} \subset c_1\sqrt{N}B_2^N$, Theorem \ref{thm:decomposition} implies that for large $N$,
$$
P_\sigma F \subset c_2(\delta)d_{\psi_\alpha}\sqrt{N}B_2^N.
$$
Hence, if we are in a situation where the $\psi_\alpha$ and the $L_2$ metrics are equivalent, then Theorem \ref{thm:decomposition} is stronger than Theorem \ref{thm:Rud}, since the $\sqrt{N}B_1^N$ component is not needed.

In fact, the gap between the two results can be considerable. As an example, let $F=\{\inr{x,\cdot} : x \in S^{n-1} \}$ and set $\mu$ to be the canonical gaussian measure on $\R^n$. Then, by Theorem \ref{thm:decomposition}, with high probability
$$
P_\sigma F \subset c(\delta) \left(\sqrt{n} B_2^N + \lambda B_\infty^N \cap c B_{\psi_2^N} \right).
$$
On the other hand, since
$$
\frac{1}{\sqrt{N}} \E \sup_{f \in F} | \sum_{i=1}^N \eps_i f(X_i) | \sim \frac{1}{\sqrt{N}} \E (\sum_{i=1}^N \|X_i\|^2_{\ell_2^n} )^{1/2} \sim \sqrt{n},
$$
then Theorem \ref{thm:Rud} only yields that
$$
P_\sigma F \subset c(\delta) \left(\sqrt{nN}B_1^N+\sqrt{N}B_2^N\right),
$$
which is a much weaker estimate.

The reason for the gap between the results is that Theorem \ref{thm:decomposition} is tailored for situations in which one has additional information on the tails of functions in the class, and in return gets more structural information
on the peaky part of coordinate projections. On the other hand, the assumptions of Theorem \ref{thm:Rud} can only give little information on the peaky part of $F$, which is captured by the $\ell_1^N$ component of that decomposition. Indeed, at best, for a fixed, reasonable class of functions $F$, one may expect that
$$
\frac{1}{\sqrt N} \E \sup_{f \in F} \left|\sum_{i=1}^N \eps_i f(X_i) \right| \leq c(F).
$$
Thus Theorem \ref{thm:Rud} only yields
$$
P_\sigma F \subset c(\delta) \left(c(F) \sqrt{N}B_1^N + d_{L_2}\sqrt{N}B_2^N\right),
$$
but with no further information on the way the coordinates are distributed in the $\ell_1^N$ component of the decomposition.
In the geometric applications we are interested in, the constant $c(F)$ grows with the ``dimension" of the class (in the example presented above, $c(F) \sim \sqrt{n}$), making Theorem \ref{thm:Rud} too weak for the analysis of such problems.

We end this section with a formulation of a simple application of the proof of Theorem
\ref{thm:decomposition}.

\begin{Corollary} \label{cor:bounded concentration}
There exist absolute constants $c_0$, $c_1$, $c_2$ and $c_3$ for
which the following holds. Let $F$ be a class of mean-zero functions
and for every $N$ set $\lambda = c_0d_{\psi_1} \max\{ \log
(c_0d_{\psi_1}N^{1/2}/\gamma_2(F,\psi_2),1)\}$. Then, for every $t
\geq c_1$, with probability at least $1-2\exp(-c_2\min \{t \log N,
t^2\})$,
$$
\sup_{f \in F} \left|\frac{1}{N} \sum_{i=1}^N f^2(X_i) - \E f^2
\right| \leq c_3  t^2\max \left\{ \lambda
\frac{\gamma_2(F,\psi_2)}{\sqrt{N}},
\frac{\gamma_2^2(F,\psi_2)}{N} \right\}.
$$
\end{Corollary}

It is important to note that using the $L_\infty$ bound to obtain a
concentration result for $F_2$ (as one does in Corollary
\ref{cor:bounded concentration}) leads to a logarithmic
looseness. Indeed, to obtain the correct estimate on the expectation
of $\sup_{f \in F} |P_N f^2 - P f^2|$ one has to truncate functions at a
level $\sim d_{\psi_1}$. This is impossible even if one considers a
single gaussian random variable. It is true that for small values of
$N$ -- when $d_{\psi_1} N^{1/2} \ll \gamma_2(F,\psi_2)$, the level
of truncation is the required one, but the resulting estimate on
$\sup_{f \in F} |P_N f^2 - P f^2|$ is trivial. Indeed, for those
values of $N$ there is no real concentration and the bound reflects
an estimate on the empirical diameter $\sup_{f \in F} (P_N
f^2)^{1/2}$. On the other hand, when $d_{\psi_1} N^{1/2} \sim
\gamma_2(F,\psi_2)$ and beyond, one starts seeing true
concentration, but then the best possible level of truncation for
those values of $N$ is off by a logarithmic factor from the required
one. Thus, even with a sharp decomposition theorem at our disposal,
a contraction based estimate on the empirical process indexed by
$F^2$ leads to a superfluous $\log N$ factor. Despite that, this
type of a decomposition argument is strong enough for many
applications (see, for example, \cite{Bour,GiaMil,Men:weak}, and
most notably, in \cite{ALPT}), because in those cases the all the
required information is when $d_{\psi_1} \sqrt{N}$ is proportional to
the complexity parameter of the class, rather than for larger values
of $N$.

If one wishes to obtain the correct estimate on $\sup_{f \in F} |P_N
f^2 - P f^2|$ for larger values of $N$, more accurate information on
the ``bounded part" of $F$ is needed. This is not surprising because
decomposition theorems like Theorem \ref{thm:decomposition} are
based solely on deviation estimates and on bounds on the $\ell_2$
norms of monotone rearrangements of $(f(X_i))_{i=1}^N$. On the other
hand, the correct rates require some sort of ``local" concentration
bounds, and those are at the heart of the proof of Theorem A.

\section{From a bounded diameter to concentration} \label{sec:conc}

Here we will remove the superfluous logarithmic factor and prove
Theorem A, by showing if $F$ is a symmetric class of mean-zero
functions, then with high probability and in expectation
\begin{equation} \label{eq:p-power-est}
\sup_{f \in F} \left|\frac{1}{N} \sum_{i=1}^N f^2(X_i) - \E f^2
\right| \leq c \max\left\{d_{\psi_1}
\frac{\gamma_2(F,\psi_2)}{\sqrt{N}}, \frac{\gamma_2^2(F,\psi_2)}{N}
\right\}.
\end{equation}
In particular, in the non-trivial range where there is actual
concentration, the dominating term is $d_{\psi_1}
\gamma_2(F,\psi_2)/{\sqrt{N}}$, which is a contraction type estimate
with the maximal norm in $\psi_1$ taking the role of the maximal
norm in $L_\infty$.

The source of difficulty in the proof of Theorem A is that the
desired concentration does not follow from the individual
concentration of each $N^{-1}\sum_{i=1}^N f^2(X_i)$ around its mean.
Rather, it is a combination of two components. First, a tail estimate on the diameter
of the ``ends" of chains
$$
\left(\sup_{f \in F} \frac{1}{N} \sum_{i=1}^N (f-\pi_{\tau_N}f)^2(X_i) \right)^{1/2},
$$
whose role
in the chaining process is to capture the
``peaky behavior" of $F$ that prevents concentration.
The second component is an analysis of the Bernoulli process
$\sup_{f \in F} \sum_{i=1}^N \eps_i (\pi_{\tau_N}f)^2(X_i)$,
conditioned on $(X_i)_{i=1}^N$. It captures the part of $F$ in which there is concentration.
Moreover, the analysis of both parts has to be carried out without resorting to a ``global" contraction argument, because
the $L_\infty$ or $\psi_2$ diameters of the relevant sets may be too large.

As a starting point of the proof of Theorem A, consider an almost optimal admissible sequence of $F$ with respect to the $\psi_2$ metric. Let $\tau_N$ be the integer $s$ satisfying that $N/2 \leq 2^s < N$. One can show \cite{MPT} that with high probability
$$
\sup_{f \in F} \frac{1}{N} \sum_{i=1}^N \left(f-\pi_{\tau_N}f\right)^2(X_i) \lesssim \frac{\gamma_2^2(F,\psi_2)}{N},
$$
which is of the desired order of magnitude. This estimate is based on Bernstein's inequality, which implies that $(N^{-1}\sum_{i=1}^N h^2(X_i))^{1/2}$ behaves like a sum of i.i.d. $\psi_2$ random variables for ``large" deviations.

Next, one has to study
$$
\sup_{f \in F} \left|\frac{1}{N} \sum_{i=1}^N (\pi_{\tau_N}f)^2(X_i) - \E (\pi_{\tau_N}f)^2 \right |,
$$
which, by a symmetrization argument behaves like
$$
\sup_{f \in F} \left|\frac{1}{N} \sum_{i=1}^N \eps_i\left((\pi_{\tau_N}f)^2(X_i) -(\pi_0 f)^2(X_i)\right) \right|.
$$
To analyze this Bernoulli process one uses a chaining argument with the same, non-random admissible sequence,
and thus one has to study the increments
\begin{equation} \label{eq:local-increment}
Pr_\eps \left(\left| \sum_{i=1}^N \eps_i \left((\pi_s f)^2 - (\pi_{s-1} f)^2 \right)(X_i) \right| > t \right),
\end{equation}
conditioned on $(X_i)_{i=1}^N$.
At every level $s \leq \tau_N$ one has to control the $2^{2^{s+1}}$ vectors in $\R^N$ of the form $(y_i)_{i=1}^N=((\pi_s f-\pi_{s-1}f)(\pi_s f+\pi_{s-1}f)(X_i))_{i=1}^N$.

Since $\pi_s f\in F$ and thanks to Theorem B, one has very accurate information on the coordinate structure of $((\pi_s f+\pi_{s-1} f)(X_i))_{i=1}^N$. However, a similar result is required for the differences $((\pi_s f-\pi_{s-1}f)(X_i))_{i=1}^N$ which takes into account the $\psi_2$ distance between $\pi_s f$ and $\pi_{s-1} f$. The desired estimate is proved in Lemma \ref{lemma-local-main} below.

Finally, to bound \eqref{eq:local-increment}, observe that for every $\ell \leq N$ and every $(y_i)_{i=1}^N$,
$$
Pr \left( \left| \sum_{i=1}^N \eps_i y_i \right| \geq \sum_{i=1}^\ell y_i^* + \sqrt{t} \left(\sum_{i=\ell+1}^N (y^2)_i^*\right)^{1/2} \right) \leq 2\exp(-t^2/2).
$$
This observation is used at the level $s$ of the chaining process for $t \sim 2^s$ and for different values of $\ell$ that depend both on $s$ and on the structure of each $(y_i)_{i=1}^N=((\pi_s f-\pi_{s-1}f)(\pi_s f+\pi_{s-1}f)(X_i))_{i=1}^N$. The crucial point in determining $\ell$ is the number of coordinates on which $\left((\pi_s f -\pi_{s-1} f)(X_i)\right)_{i=1}^N$ does not ``behaves regularly" in the sense of Theorem \ref{thm:decomposition}.

We begin the proof with a ``local" version of Theorem \ref{thm:main-est} -- for a finite class $H$, in
which for every $h \in H$ the bound on $(\sum_{i \in I} h^2(X_i))^{1/2}$ is given using
$\|h\|_{\psi_2}\sqrt{\log |H|}$ and $\|h\|_{\psi_1}$ rather than using the global
parameters $\gamma_2(H,\psi_2)$ and $d_{\psi_1}$ that are used in Theorem B.

\begin{Lemma} \label{lemma-local-main}.
There exists absolute constants $c_1, c_2$, $c_3$ and $c_4$ for
which the following holds. Let $H$ be a class of mean-zero functions
and set $k=\log |H|$. Then, for every $u \geq c_1$, with probability
at least $1-2\exp(-c_2 \max\{k,\log N\} u)$, for every $I \subset
\{1,...,N\}$ and every $h \in H$,
$$
\left(\sum_{i \in I} h^2(X_i)\right)^{1/2} \leq c_4 u \left(
\|h\|_{\psi_2}\sqrt{k} + \|h\|_{\psi_1} \sqrt{|I|} \log (eN/|I|)
\right).
$$
An analogous result holds for any $\psi_\alpha$ norm for $1 <
\alpha \leq 2$, with $\log^{1/\alpha}(eN/|I|)$ taking the place of
$\log(eN/|I|)$.
\end{Lemma}
We will prove the lemma for $\alpha=1$ since this is the only case
we will actually use. The proof for $1 < \alpha \leq 2$ follows the
same lines and is omitted.

The proof of Lemma \ref{lemma-local-main} is very similar in nature
to the proof of Theorem \ref{thm:main-est} and will use its
notation. Again, we will denote by $E_m$ the collection of subsets
of $\{1,...,N\}$ of cardinality $m$.

\proof Recall that for every $h \in H$,
$$ \sup_{|I|=m}
\left(\sum_{i \in I} h^2(X_i) \right)^{1/2} \leq C \sup_{v \in B_m}
\sum_{i=1}^N v_i h(X_i)$$ where $B_m$ was defined in \eqref{eq:B_m}.
First, assume that $|H| \geq N$ and that $m=2^{r_0}$ satisfies that
$\kappa_0m \log(eN/m) \geq \max\{ \log |E_m|, \ \log |B_m|\}$. We
can assume without loss of generality that $\log |E_m| \gtrsim \log
|H|$. Indeed, if $\log |E_m| \lesssim \log |H|$ then the required
estimate follows easily from a $\psi_2$ estimate and the union
bound, since $\log (|H| \cdot |B_m|) \leq c_0 \log |H|$.

Recall that for every $v
\in B_m$, $|{\rm supp}(v)| \leq m$ and that there is a set $I_{\ell_1}$
of cardinality $m/2$ such that $\|P_{I_{\ell_1}}v\|_\infty \leq
1/(m/2)^{1/2}$. Let $J_1$ be the complement of $I_{\ell_1}$ in
$|{\rm supp}(v)|$, and so on for $\ell_r=m/2^r$, $r \leq r_1$, where
$r_1$ will be named later. Observe that for every $v \in B_m$,
$P_{J_r} v \in B_{\ell_{r}}$. Since $\max_{i \in I_{\ell_r}} \|v_i
h(X_i)\|_{\psi_1} \leq \|h\|_{\psi_1}/\sqrt{\ell_r}$, then by
Bernstein's inequality, for every $u_1$ larger than an absolute
constant,
\begin{align} \label{eq:(*)}
& Pr  \left( \exists h \in H, \ v \in B_m, \ r \leq r_1  \ \ \left|
\sum_{i \in I_{\ell_r}} v_i h(X_i) \right| \geq
u_1\|h\|_{\psi_1}\sqrt{\ell_{r}}\log(eN/\ell_{r}) \right)
\nonumber
\\
& \leq  2|H| \sum_{r=1}^{r_1} |B_{\ell_{r}}| \exp(-c_1 u_1
\ell_{r} \log(eN/\ell_{r})) \nonumber
\\
& \leq 2|H| \exp(-c_2 u_1 \ell_{r_1} \log(eN/\ell_{r_1}))=(*).
\end{align}
Since $\ell_r=m/2^r$, we set $r_1$ to be the largest integer for
which
$$ ({m}/{2^{r_1}}) \log (eN/(m/2^{r_1})) \gtrsim \log |H|$$ and
since $\log |E_m| \gtrsim \log |H|$ such an integer exists. Thus,
for $u_1 \geq c_3$ it is evident that $ (*) \leq 2\exp(-c_4u_1 \log
|H|)$.

Next, for every $v \in B_m$ consider the projection $P_{J_{r_1}}v$.
Since $\|P_{J_{r_1}}v\|_2 \leq 1$ then $\|\sum_{i \in J_{r_1}} v_i
h(X_i) \|_{\psi_2} \leq c_5\|h\|_{\psi_2}$. Therefore,
\begin{align*}
Pr & \left( \exists \ v \in B_m, \  h \in H \ : \ \left| \sum_{i \in
J_{r_1}} v_i h(X_i) \right| \geq u_2 \sqrt{\log |H|} \|h\|_{\psi_2}
\right)
\\
\leq & 2|H| \cdot |B_{m/2^{r_1}}| \exp(-c_6 u_2^2 \log |H|) \leq
\exp(-c_7 u^2_2 \log |H|),
\end{align*}
provided that $u_2 \geq c_8$.

Therefore, if $u$ is sufficiently large, then with probability at
least $1-2\exp(-c_9 u \log |H|)$, for every $h \in H$
\begin{align*}
& \sup_{v \in B_m}  \sum_{i=1}^N v_i h(X_i)
\\
\leq & c_{10}u  \left(\|h\|_{\psi_2} \sqrt{\log{H}} +
\|h\|_{\psi_1}\sum_{r=1}^{r_1} \sqrt{\ell_{r}}\log(eN/\ell_{r})
\right)
\\
\leq & c_{11}u \left(\|h\|_{\psi_2} \sqrt{\log{H}} +
\|h\|_{\psi_1}\sqrt{m}\log(eN/m) \right).
\end{align*}
Since $1 \leq m \leq N$ and $|H| \geq N$, the claim holds for
any such $m$.

Now, assume that $|H| < N$. Then, set $r_1=r_0$, and by
\eqref{eq:(*)}, with probability at least $1-2\exp(-c_{12}u\log N)$,
for every $h \in H$ and $v \in B_m$, $\left| \sum_{i=1}^N v_i h(X_i)
\right| \leq c_{13}u\|h\|_{\psi_1} \sqrt{m}\log(eN/m)$. Again,
summing the probabilities over every $1 \leq m \leq N$ the claim
follows.
\endproof

Recall that $\tau_N$ is the integer $s$ for which $N/2 \leq 2^s < N$. The
sets $H$ we will be interested in are the sets of links at the level
$s$, namely, $\Delta_s=\{\Delta_s(f)=\pi_sf-\pi_{s-1}f : f \in F\}$,
for $s \leq \tau_N$, where $(F_s)_{s \geq 0}$ is
an almost optimal admissible sequence of $F$ with respect to the
$\psi_2$ norm.

Let us summarize the information we have on the set
$$
\left\{ \left( (\Delta_s(f))(X_i)\cdot(\pi_s f +\pi_{s-1} f)(X_i)\right)_{i=1}^N
: f \in F \right\}.
$$
Consider the following events: let
\begin{align} \label{eq:A_t}
& A_t = \Biggl\{ (X_i)_{i=1}^N : \forall I \subset \{1,...,N\}, \
\\
& \sup_{f \in F} \left(\sum_{i \in I} f^2(X_i) \right)^{1/2} \leq
\kappa_1 t\left(\gamma_2(F,\psi_2) + d_{\psi_1}
\sqrt{|I|}\log(eN/|I|) \right) \Biggl\}, \nonumber
\end{align}
and
\begin{align} \label{eq:B_t^s}
& B_t^s= \Biggl\{  (X_i)_{i=1}^N :  \forall f \in F, \ \ \forall I
\subset \{1,...,N\},
\\
& \left(\sum_{i \in I} (\Delta_s f)^2(X_i) \right)^{1/2} \leq
\kappa_1 t\left(2^{s/2}\|\Delta_s(f)\|_{\psi_2} +
\|\Delta_s(f)\|_{\psi_1} \sqrt{|I|}\log(eN/|I|) \right) \Biggl\},
\nonumber
\end{align}
where $\kappa_1$ is a suitable absolute constant.

By Theorem \ref{thm:main-est}, for every $t \geq c_{1}$, $Pr(A_t)
\geq 1-2\exp(-c_{2} t \log N)$, while applying Lemma
\ref{lemma-local-main} it is evident that for every $t \geq c_3$ and
every $s \leq \tau_N$, $Pr(B_t^s) \geq 1-2\exp(-c_4\max\{2^s,\log
N\}t)$.

Let us consider the Bernoulli process $\left|\sum_{i=1}^N \eps_i
((\pi_{\tau_N}f)^2(X_i)-f_0^2(X_i)) \right|$ conditioned on the set
$\Omega_t = A_t \cap \left(\bigcap_{s \leq \tau_N} B_t^s \right)$. Observe that on $\Omega_t$ we have enough information to identify the cardinality of each set of ``large" coordinates of individual functions $\Delta_s(f)$, and thus the point from which each vector $((\Delta_s(f))(X_i))_{i=1}^N$ behaves regularly. The parameter we will use to identify the point from which the regular behavior begins is
$$
m\left(\Delta_s(f)\right)=\min\left\{ m : 2^{s/2}\|\Delta_s(f)\|_{\psi_2} \leq \|\Delta_s(f)\|_{\psi_1} \sqrt{m} \log(eN/m) \right\}.
$$

\begin{Theorem} \label{thm:local-Bernoulli}
There exist absolute constants $c_1,c_2,c_3$ and $c_4$ for which the
following holds. If $f_0 \in F$ then for every $t \geq c_1$, $u \geq
c_2$ and $(X_i)_{i=1}^N \in \Omega_t$,
\begin{equation*}
Pr_\eps \left(\exists f \in F,  \ \ \left|\sum_{i=1}^N \eps_i \left(
(\pi_{\tau_N}f)^2(X_i) -f_0^2(X_i) \right)\right| \geq u \rho_t
\right) \leq 2\exp(-c_3u^2),
\end{equation*}
where
$$
\rho_t = c_4 t^2 \left(\sqrt{N}d_{\psi_1}\gamma_2(F,\psi_2)+
\gamma_2^2(F,\psi_2)\right).
$$
\end{Theorem}
For the proof we will need the following definition.
\begin{Definition} \label{def:dominating}
Let $u=(u_i)_{i=1}^N$, $I \subset \{1,...,N\}$  and
$v=(v_i)_{i=1}^{|I|}$. We say that $v$ dominates $u$ on $I$ if for
every $i \in I$, $(P_I u)_i^* \leq v_i^*$. In other words, if a
monotone rearrangement of $P_I u$ is smaller than that of $v$
coordinate-wise on $I$.
\end{Definition}

\proof Fix $f_0 \in F$ and for every $f \in F_{\tau_N}$ write
$f^2-f_0^2=\sum_{s=1}^{\tau_N} (\pi_s f)^2 - (\pi_{s-1} f)^2$. Let
$1 \leq s \leq \tau_N$ and consider a link $(\pi_s f)^2 - (\pi_{s-1}
f)^2$. Set $h_-=\pi_{s} f-\pi_{s-1}f$,
$h_+=\max\{\pi_{s}f,\pi_{s-1}f\}$ and for every $(X_i)_{i=1}^N$ let
$v_-=(h_-(X_i))_{i=1}^N$ and $v_+=(h_+(X_i))_{i=1}^N$. Also, for
every $h_-$, recall that $m(h_-)$ is the smallest integer such that
$2^{s/2}\|h_-\|_{\psi_2} \leq \|h_-\|_{\psi_1} \sqrt{m}\log(eN/m)$
and if the smallest is $m \geq N$, set $m(h_-)=N$.

Since $(X_i)_{i=1}^N \in B_t^s$ then for every $I
\subset \{1,...,N\}$
$$
\left(\sum_{i \in I} h_-^2(X_i) \right)^{1/2} \leq \kappa_1t
\left(2^{s/2}\|h_-\|_{\psi_2} + \|h_-\|_{\psi_1}
\sqrt{|I|}\log(eN/|I|) \right),
$$
and let us consider two cases. The first is when $m(h_-) \leq 2^s $
and the second is when the reverse inequality holds.

To handle the
first case, when $m(h_-) \leq 2^s$, observe that by the subgaussian
inequality for Bernoulli sums, for every $u>0$, with probability at
least $1-2\exp(-c_1u^22^s)$
\begin{align*}
\left| \sum_{i=1}^N \eps_i ((\pi_s f)^2 - (\pi_{s-1} f)^2 )\right| &
\leq \sum_{i=1}^{2^s} (v_-v_+)_i^* + |\sum_{i > 2^s} \eps_i
(v_-v_+)_i^* |
\\
& \leq \sum_{i=1}^{2^s} (v_-v_+)^*_i + u2^{s/2}\left(\sum_{i
> 2^s} (v_-^2v_+^2)^*_i \right)^{1/2}
\end{align*}
where, as always, $(x^*_i)_{i \geq 1}$ denotes a non-increasing
rearrangement of $(|x_i|)_{i \geq 1}$.

Clearly, $ \sum_{i=1}^{2^s} (v_-v_+)^*_i \leq \left(\sum_{i=1}^{
2^s} (v_-^2)^*_i\right)^{1/2} \left(\sum_{i=1}^{2^s}
(v_+^2)^*_i\right)^{1/2}. $ Since $(v_+)_i \leq \max\{|(\pi_s
f)(X_i)|,|(\pi_{s-1} f)(X_i)|\}$ and $(X_i)_{i=1}^N \in A_t$ then
\begin{align*}
\left(\sum_{i=1}^{2^s} (v_+^2)^*_i\right)^{1/2} \leq & 2\sup_{f
\in F} \sup_{|I|=2^s} \left(\sum_{i \in I} f^2(X_i) \right)^{1/2}
\\
\leq & c_2t\left(\gamma_2(F,\psi_2) + d_{\psi_1}2^{s/2}\log (eN/2^s)
\right).
\end{align*}
Also, because $m(h_-) \leq 2^s$ and $(X_i)_{i=1}^N \in B_t^s$, it is
evident that
\begin{align*}
\left(\sum_{i=1}^{2^s} (v_-^2)^*_i\right)^{1/2} \leq &
t\kappa_1\left(2^{s/2}\|h_-\|_{\psi_2} + \|h_-\|_{\psi_1}
2^{s/2}\log(eN/2^s) \right)
\\
\leq &  c_3t \|h_-\|_{\psi_1} 2^{s/2}\log(eN/2^s).
\end{align*}
Hence, recalling that
$$
2^{s/2}\|h_-\|_{\psi_1} = 2^{s/2}\|\Delta_s
(f)\|_{\psi_1} \leq 2^{s/2}\|\Delta_s (f)\|_{\psi_2} \leq
\gamma_2(F,\psi_2),
$$
one has
\begin{align} \label{eq:small-m}
& \sum_{i=1}^{2^s} (v_-v_+)^*_i \nonumber
\\
\leq & c_4 t^2 \left( \gamma_2(F,\psi_2) \|h_-\|_{\psi_1}
2^{s/2}\log(eN/2^s) + d_{\psi_1}\|h_-\|_{\psi_1} 2^{s}\log^2(eN/2^s)
\right) \nonumber
\\
\leq & 2c_4t^2d_{\psi_1} \gamma_2(F,\psi_2) 2^{s/2}\log^2(eN/2^s).
\end{align}

Next, let us consider the term
$$
\left(\sum_{i \geq 2^s} (v_-^2v_+^2)^*_i \right)^{1/2} =
\inf_{\{J:|J|=N-2^s\}} \left(\sum_{i \in J} (v_-^2v_+^2)^*_i
\right)^{1/2}.
$$
Let $J$ be the set of the $N-2^s$ smallest coordinates of
$v_-$. Observe that for every $I \subset \{1,...,N\}$, $|I| \geq 2^s$ one has
\begin{align*}
\left(\sum_{i \in I} (v_-^2)_i \right)^{1/2} & \leq \kappa_1 t
\left(2^{s/2}\|h_-\|_{\psi_2}+ \|h_-\|_{\psi_1} \sqrt{|I|}
\log(eN/|I|)\right)
\\
& \leq 2\kappa_1t\|h_-\|_{\psi_1} \sqrt{|I|} \log(eN/|I|),
\end{align*}
since $m(h_-) \leq 2^s$.

Thus, for every $i > 2^s$,
$$
(v_-)_i^* \leq 2\kappa_1t
\|h_-\|_{\psi_1} \log(eN/i),
$$
and in particular, $v_-$ is dominated by
$(2\kappa_1 t \|h_-\|_{\psi_1} \log(eN/i))_{i > 2^s}$ on the
set $J$ and
$$
\|P_J v_-\|_\infty \leq 2\kappa_1t \|h_-\|_{\psi_1} \log(eN/2^s).
$$
To obtain a similar control over the vector $v_+$, let $m_0$ be the smallest integer such that $\gamma_2(F,\psi_2) \leq
d_{\psi_1} \sqrt{m} \log(eN/m)$, and if the smallest one is larger
than $N$, set $m_0=N$. Just as we did for $v_-$, if $f \in F$ and
$(X_i)_{i=1}^N \in A_t$, and if we set
$(u_i)_{i=1}^N=(f(X_i))_{i=1}^N$, then for every $i \geq m_0$,
$u^*_i \leq 2t\kappa_1d_{\psi_1} \log(eN/i)$. Therefore, if $I_+$ is
the set of the $m_0$ largest coordinates of $(v_+)_{i=1}^N$, then
$(u_i)_{i=1}^N$ is dominated by $(2t\kappa_1d_{\psi_1}\log(eN/i))_{i
> m_0}$ on $I_+^c$. Therefore,
\begin{align*}
\left(\sum_{i \in J \cap I_+^c} (v_-^2v_+^2)_i\right)^{1/2} \leq &
8\kappa_1^2t^2 d_{\psi_1} \|h_-\|_{\psi_1} \left(\sum_{i=1}^N
\log^4(eN/i)\right)^{1/2}
\\
\leq & c_5 t^2 d_{\psi_1} \|h_-\|_{\psi_1} \sqrt{N},
\end{align*}
implying that

\begin{align*}
& \left( \sum_{i > 2^s} (v_-^2v_+^2)^*_i\right)^{1/2} \leq
\left(\sum_{i \in J \cap I_+} (v_-^2v_+^2)_i\right)^{1/2}
+\left(\sum_{i \in J \cap I_+^c} (v_-^2v_+^2)_i\right)^{1/2}
\\
\leq & \left(\sum_{i=1}^{m_0} (v_+^*)_i^2 \right)^{1/2} \|P_J
v_-\|_\infty + c_5t^2 \|h_-\|_{\psi_1} d_{\psi_1} \sqrt{N}
\\
\leq & c_{6}t^2 \left(\gamma_2(F,\psi_2) \|h_-\|_{\psi_1}
\log(eN/2^s) + d_{\psi_1} \|h_-\|_{\psi_1} \sqrt{N} \right).
\end{align*}
Thus, if $m(h_-) \leq 2^s$ then

\begin{align} \label{eq:first-estimate}
\sum_{i=1}^{2^s} & (v_-v_+)_i^* + u 2^{s/2} \left(\sum_{i \geq
2^s} (v_-^2v_+^2)^*_i \right)^{1/2} \nonumber
\\
\leq & c_{7}  t^2 2^{s/2}\Biggl( d_{\psi_1} \gamma_2(F,\psi_2)
\log^2(eN/2^s) \nonumber
\\
+ & u \left(\gamma_2(F,\psi_2) \|h_-\|_{\psi_1}
\log(eN/2^s)+d_{\psi_1} \|h_-\|_{\psi_1} \sqrt{N}\right) \Biggr)
\nonumber
\\
\leq & c_{8} u t^2 d_{\psi_1} \left(  \gamma_2(F,\psi_2) 2^{s/2}
\log^2(eN/2^s) + 2^{s/2} \|h_-\|_{\psi_1}\sqrt{N} \right),
\end{align}
provided that $u \geq 1$.

Next, we turn to the case when $m(h_-) > 2^s$. Let $I_-$ be the
set of the $m(h_-)$ largest coordinates of $v_-$ and again, $I_+$
is the set of the $m_0$ largest coordinates of $v_+$. Therefore,
if $I \subset \{1,...,N\}$ and $|I| \geq m(h_-)$ then
$$
\left(\sum_{i \in I} (v_-^2)_i \right)^{1/2} \leq 2\kappa_1t
\|h_-\|_{\psi_1} \sqrt{|I|} \log(eN/|I|),
$$
and thus, for $i \geq m(h_-)$,
$$
(v_-)^*_i \leq 2\kappa_1t
\|h_-\|_{\psi_1} \log(eN/i).
$$
Also, from the definition of $m(h_-)$
it is evident that
$$
\left(\sum_{i=1}^{m(h_-)} (v_-^2)_i^* \right)^{1/2} \leq
4\kappa_1t 2^{s/2} \|h_-\|_{\psi_2},
$$
implying that
$$
\|P_{I_-^c}v_-\|_\infty \leq 4\kappa_1t (2^s/m(h_-))^{1/2}
\|h_-\|_{\psi_2} \leq c_9t\|h_-\|_{\psi_2},
$$
because $2^s \leq
m(h_-)$.

Set $I=I_- \cup I_+ = (I_+ \backslash I_-) \cup I_-$. Note that
\begin{align*}
\sum_{i \in I_-} |(v_-v_+)_i| & \leq  \left(\sum_{i \in I_-}
(v_-^2)_i\right)^{1/2}\left(\sum_{i \in I_-} (v_+^2)_i\right)^{1/2}
\\
& \leq c_{10}t^22^{s/2}\|h_-\|_{\psi_2}
\left(\gamma_2(F,\psi_2)+d_{\psi_1}\sqrt{N}\right),
\end{align*}
where we have used that
\begin{align*}
\left(\sum_{i \in I_-} (v_+^2)_i \right)^{1/2} & \leq 2\sup_{f \in F} \left(\sum_{i=1}^N f^2(X_i)\right)^{1/2}
\\
& \lesssim
t\left(\gamma_2(F,\psi_2)+d_{\psi_1}\sqrt{N}\right).
\end{align*}

Moreover, applying the bound on $\|P_{I_-^c}v_-\|_{\infty}$,
\begin{align*}
& \left(\sum_{I_+ \backslash I_-} (v_-^2v_+^2)_i \right)^{1/2} \leq
\|P_{I_-^c} v_-\|_\infty \left(\sum_{i=1}^N (v_+^2)_i \right)^{1/2}
\\
& \leq c_{11} t^2 \|h_-\|_{\psi_2}
\left(\gamma_2(F,\psi_2)+d_{\psi_1}\sqrt{N}\right).
\end{align*}
That leaves us with the coordinates that are outside $I$, that is,
outside both $I_-$ and $I_+$. Observe that $v_-$ is dominated on
$I^c$ by $(2\kappa_1t\|h_-\|_{\psi_1} \log(eN/i))_{i=1}^{|I^c|}$
and $v_+$ is dominated on $I^c$ by
$(2\kappa_1td_{\psi_1}\log(eN/i))_{i=1}^{|I^c|}$. Hence,
\begin{align*}
\left(\sum_{i \in I^c} (v_-^2v_+^2)_i \right)^{1/2} & \leq
4\kappa_1^2t^2\|h_-\|_{\psi_1} d_{\psi_1} \left(\sum_{i=1}^N
\log^4(eN/i) \right)^{1/2}
\\
& \leq c_{11} t^2\|h_-\|_{\psi_1} d_{\psi_1} \sqrt{N}.
\end{align*}
Therefore, if $m(h_-) \geq 2^s$, then with
$(\eps_i)_{i=1}^N$-probability at least $1-2\exp(-c_{12}u^22^s)$,
\begin{align} \label{eq:second-estimate}
\left| \sum_{i=1}^N \eps_i (v_-v_+)_i \right| & \leq c_{13}ut^2
2^{s/2} \Biggl(\|h_-\|_{\psi_2}
\left(\gamma_2(F,\psi_2)+d_{\psi_1}\sqrt{N}\right) +
\|h_-\|_{\psi_1}d_{\psi_1} \sqrt{N} \Biggr) \nonumber
\\
& \leq c_{14}u t^2 2^{s/2} \|h_-\|_{\psi_2}
\left(\gamma_2(F,\psi_2)+d_{\psi_1}\sqrt{N}\right),
\end{align}
provided that $u \geq 1$.

Combining \eqref{eq:first-estimate} and \eqref{eq:second-estimate},
and since there are at most $2^{2^{s+1}}$ links at the $s$-level, it
is evident that for every $t,u \geq c_{15}$ and every $s \leq
\tau_N$,
\begin{align} \label{eq:s-estimate}
Pr_\eps & \left( \exists f \in F : \left|\sum_{i=1}^N \eps_i \left(
(\pi_s(f))^2(X_i) - (\pi_{s-1}(f))^2(X_i) \right) \right| \geq
\rho(s,f) u \Big| \Omega_t \right) \nonumber
\\
\leq & 2\exp(-c_{16}u^2 2^s),
\end{align}
where
\begin{align*}
\rho(s,f) & \sim u t^2 d_{\psi_1}  \left(  \gamma_2(F,\psi_2)
2^{s/2}\log^2(eN/2^s) + 2^{s/2}\|\Delta_s(f)\|_{\psi_1}\sqrt{N}
\right)
\\
& + ut^2 \left(2^{s/2}\|\Delta_s(f)\|_{\psi_2}\right)
\left(\gamma_2(F,\psi_2)+d_{\psi_1}\sqrt{N}\right).
\end{align*}
It remains to show that for every $f \in F$,
$$
\sum_{\{s : 2^s \leq N\}} \rho(s,f) \lesssim ut^2
\left(\sqrt{N}d_{\psi_1}\gamma_2(F,\psi_2)+
\gamma_2^2(F,\psi_2)\right),
$$
which is straightforward because for an almost optimal admissible
sequence,
$$
\sum_{s \geq 1} 2^{s/2} \|\Delta_s(f)\|_{\psi_1} \leq
2\gamma_2(F,\psi_2) \ \ {\rm and} \ \ \sum_{\{s : s \leq \tau_N\}}
2^{s/2} \log^2(eN/2^s) \sim \sqrt{N}.
$$
\endproof

We need an additional preliminary result which allows one to move freely between the empirical process and the Bernoulli one -- the Gin\'{e}-Zinn
symmetrization Theorem \cite{GZ84}:
\begin{Theorem} \label{thm:GZ}
Let $H$ be a class of functions and set $\alpha^2=\sup_{h \in H}
\E(h-\E h)^2$. For every integer $N$ and any $t \geq 2^{1/2}\alpha
N^{1/2}$,
\begin{equation*}
Pr \left( \sup_{h \in H}  \left|\sum_{i=1}^N (h(X_i)-\E h) \right|
>t \right)
\leq 4 Pr\left( \sup_{h \in H} \left|\sum_{i=1}^N \eps_i h(X_i)
\right| > t/4 \right).
\end{equation*}
\end{Theorem}

Combining Theorem \ref{thm:GZ} and Theorem \ref{thm:local-Bernoulli}
we obtain the next result on the ``beginning" of every chain.
\begin{Theorem} \label{thm:main-concentration-mid-part}
There exist absolute constants $c_1,c_2$ and $c_3$ for which the
following holds. Let $F$ be a class of mean-zero functions and let
$(F_s)_{s \geq 0}$ be an almost optimal admissible sequence with
respect to the $\psi_2$ norm. Then, for every $f_0 \in F$ and $x
\geq c_1$, with probability at least $1-2\exp(-c_2x^{2/5})$,
\begin{align*}
\sup_{f \in F} & \left| \sum_{i=1}^N
(\pi_{\tau_N}f)^2(X_i)-f_0^2(X_i) - \E
((\pi_{\tau_N}f)^2-f_0^2(X_i)) \right|
\\
& \leq c_3x \left(\sqrt{N}d_{\psi_1}\gamma_2(F,\psi_2)+
\gamma_2^2(F,\psi_2)\right).
\end{align*}
\end{Theorem}

\begin{Remark} The power of $x^{2/5}$ in the exponent is
likely to be an artifact of the proof. We
made no effort to optimize this power since it is not of major
importance in the problems we wish to address, and because any exponential
tail estimate would give us the integrability properties we need.
\end{Remark}

\proof Fix $f_0 \in F$ and let $H=\{(\pi_{\tau_N}f)^2-f_0^2 : f \in
F\}$. It is standard to verify that  $\alpha^2 =\sup_{h \in H}
\E(h-\E h)^2 \leq c_0d_{\psi_1}^4$. Since $F$ is symmetric then
$\gamma_2(F,\psi_2) \geq d_{\psi_1}$ and $x \sqrt{N} d_{\psi_1}
\gamma_2(F,\psi_2) \geq 2 \sqrt{N} \alpha$ provided that $x \geq
2c_0$.

If we set $\rho=\left(\sqrt{N}d_{\psi_1}\gamma_2(F,\psi_2)+
\gamma_2^2(F,\psi_2)\right)$ then by Theorem \ref{thm:GZ} and the
definition of $H$, for every $x \geq 2c_0$,
\begin{align*}
Pr & \left( \sup_{f \in F}  \left|\sum_{i=1}^N
(\pi_{\tau_N}f)^2(X_i)-f_0^2(X_i) - \E ((\pi_{\tau_N}f)^2-f_0^2)
\right|
> x \rho \right)
\\
\leq & 4 Pr\left( \sup_{f \in F} \left|\sum_{i=1}^N \eps_i
((\pi_{\tau_N}f)^2(X_i)-f_0^2(X_i))\right|
> x\rho/4 \right)
\\
= & 4 \E_X Pr_\eps \left( \sup_{f \in F} \left|\sum_{i=1}^N \eps_i
((\pi_{\tau_N}f)^2(X_i)-f_0^2(X_i))\right| > x\rho/4 \right),
\end{align*}
by Fubini's Theorem.
Using the notation of \eqref{eq:A_t} and \eqref{eq:B_t^s}, for $t >
c_1$, let $\Omega_t = A_t \cap \left(\bigcap_{s \leq \tau_N} B_t^s
\right)$ and observe that
$$
Pr(\Omega_t^c) \leq Pr(A_t^c)+\sum_{s=1}^{\tau_N} Pr(B_t^c) \leq 2
\exp(-c_2t \log{N}).
$$

Thus, if we set $u_t=x/t^2$, then as long as $u_t \geq c_4$ (or in
other words, for every $t$ such that $x \geq c_4 t^2$), Theorem
\ref{thm:local-Bernoulli} implies that
\begin{align*}
& \E_X Pr_\eps \left(\sup_{f \in F} \left|\sum_{i=1}^N \eps_i
((\pi_{\tau_N}f)^2(X_i)-f_0^2(X_i))\right| \geq u_t t^2 \rho \right)
(\IND_{\Omega_t}+\IND_{\Omega_t^c})
\\
\leq & 2\left(\exp(-c_5x^2/t^4)+ \exp(-c_5t) \right) \leq 2\exp(-c_6
x^{2/5}),
\end{align*}
where the last inequality holds if we take $t=x^{2/5} \geq 1$.
\endproof

The last component in the proof of Theorem A is an estimate on the
``end" of each chain, that is, $f-\pi_{\tau_N}f = \sum_{s > \tau_N}
\Delta_s(f)$. Its proof is a combination of Bernstein's inequality
and a chaining argument (see Lemma 1.5 in \cite{MPT}), and the key
point is the observation that for every $f,g$ and every $u \geq 1$,
with probability at least $1-2\exp(-cNu^2)$, $(P_N(f-g)^2)^{1/2}
\leq u \|f-g\|_{\psi_2}$.
In particular one has
\begin{Lemma} \label{lemma:MPT} \cite{MPT}
There exist absolute constants $c_1$, $c_2$, $c_3$ and $c_4$ for
which the following holds. Let $(F_s)_{s \geq 0}$ be an almost
optimal admissible sequence of $F$ with respect to the $\psi_2$
norm. Then, for every $u \geq c_1$, with probability at least
$1-2\exp(-c_2Nu^2)$, for every $f \in F$,
$$
\sup_{f \in F} \left(P_N
(f-\pi_{\tau_N}(f))^2 \right)^{1/2} \leq c_3 u
\frac{\gamma_2(F,\psi_2)}{\sqrt{N}},
$$
and
$$
\E \sup_{f \in F}
\left(P_N(f-\pi_s(f))^2 \right)^{1/2} \leq c_4
\frac{\gamma_2(F,\psi_2)}{\sqrt{N}}.
$$
\end{Lemma}

Finally, let us reformulate Theorem A.
\begin{Theorem} \label{thm:concentration-main}
There exist absolute constants $c_1$, $c_2$, $c_3$ and $c_4$ for
which the following holds. If $F$ is a symmetric class of mean-zero
functions, then for every $x \geq c_1$, with probability at least
$1-2\exp(-c_2x^{2/5})$,
$$
\sup_{f \in F} \left|\frac{1}{N}\sum_{i=1}^N f^2(X_i) - \E f^2
\right| \leq c_3x \left(d_{\psi_1}
\frac{\gamma_2(F,\psi_2)}{\sqrt{N}} +
\frac{\gamma_2^2(F,\psi_2)}{N}\right).
$$
In particular,
$$
\E \sup_{f \in F} \left|\frac{1}{N}\sum_{i=1}^N f^2(X_i) - \E f^2
\right| \leq c_4 \left(d_{\psi_1}
\frac{\gamma_2(F,\psi_2)}{\sqrt{N}} +
\frac{\gamma_2^2(F,\psi_2)}{N}\right).
$$
\end{Theorem}
\proof Let $(F_s)_{s \geq 0}$ and $\tau_N$ be as above. Then, for
every $f \in F$,
\begin{align*}
& \sum_{i=1}^N (f^2(X_i)-\E f^2) = \sum_{i=1}^N
(f^2(X_i)-(\pi_{\tau_N}f)^2(X_i))
\\
&  +\sum_{i=1}^N \left((\pi_{\tau_N}f)^2(X_i)-\E
(\pi_{\tau_N}f)^2\right) +N \E ((\pi_{\tau_N}f)^2 - f^2)
\\
& \leq  2 \left(\sum_{i=1}^N (f-\pi_{\tau_N}f)^2(X_i) \right)^{1/2}
\sup_{f \in F} \left(\sum_{i=1}^N f^2(X_i) \right)^{1/2}
\\
& + 2N\sup_{f \in F} (\E (f-\pi_{\tau_N}f)^2)^{1/2} \cdot \sup_{f
\in F} (\E f^2)^{1/2}
\\
& + \sup_{f \in F} \sum_{i=1}^N ((\pi_{\tau_N}f)^2(X_i)-\E
(\pi_{\tau_N}f)^2).
\end{align*}
By Lemma \ref{lemma:MPT} combined with Theorem \ref{thm:main-est},
with probability at least $1-2\exp(-c_1Nt)-2\exp(-t^{1/2}\log{N})$
the first and second terms are at most
$$
c_2 t \gamma_2(F,\psi_2)
\left(\gamma_2(F,\psi_2)+d_{\psi_1}\sqrt{N}\right)
$$
for $t \geq c_3$. The third term may be bounded using Theorem
\ref{thm:main-concentration-mid-part}. Indeed, for every such $t$,
\begin{align*}
& \sum_{i=1}^N \left(((\pi_{\tau_N}f)^2(X_i)-f_0^2(X_i))-\E
(\pi_{\tau_N}f)^2-f_0^2\right)
\\
& \leq c_4t \gamma_2(F,\psi_2)
\left(\gamma_2(F,\psi_2)+d_{\psi_1}\sqrt{N}\right).
\end{align*}
with probability at least $1-2\exp(-c_5t^{2/5})$.

Finally, a similar argument to the one used in the proof of Theorem
\ref{thm:local-Bernoulli} shows that for every such $t$, with
probability at least $1-2\exp(-c_5t^{2/5})$,
\begin{align*}
\sum_{i=1}^N (f^2_0(X_i)-\E f^2_0) & \leq
c_6t(d_{\psi_2}^2+d_{\psi_1}d_{\psi_2}\sqrt{N})
\\
& \leq
c_6t(\gamma_2^2(F,\psi_2)+d_{\psi_1}\gamma_2(F,\psi_2)\sqrt{N}).
\end{align*}
Hence, with probability at least $1-4(\exp(-c_1t^{1/2}\log N
)-\exp(-c_5t^{2/5}))$,
$$
\sup_{f \in F} \left|\frac{1}{N} \sum_{i=1}^N f^2(X_i) - \E f^2
\right| \leq Ct
\left(d_{\psi_1}\frac{\gamma_2(F,\psi_2)}{\sqrt{N}}+
\frac{\gamma_2^2(F,\psi_2)}{N} \right),
$$
as required.

The claim regarding the expectation follows from an integration
argument and is omitted.
\endproof

\section{Applications} \label{sec:applications}
In this final section we will present several geometric applications
of our three main results, though as pointed out in the
introduction, there are numerous other applications in Empirical
Processes Theory, Nonparametric Statistics and Asymptotic Geometric
Analysis that will not be mentioned here.

It is well known that many results in Asymptotic Geometric Analysis
are based on a random selection argument, for example, a random
choice of a section or of a projection of a convex body in $\R^n$.
Historically, the motivation was to understand the geometry of
convex bodies and thus the models of random selection that had been
studied were rather limited. Indeed, in classical results such as
Dvoretzky's Theorem, low-$M^*$ estimates and many others (see, e.g.
\cite{MilSch,Pis:book}), the selection was preformed using a random
point on a Grassman manifold $G_{n,k}$ relative to the Haar measure,
or by applying a gaussian operator $\sum_{i=1}^k
\inr{G_i,\cdot}e_i$ to the given body, with $(G_i)_{i=1}^k$
selected independently according to the canonical gaussian measure
on $\R^n$.

In recent years, the distribution of volume in a convex body has become
a central area of interest in Asymptotic Geometric Analysis. Hence,
it is natural to ask whether the classical results in the area can
be extended to other random selection methods, endowed by these
volume measures, or, more generally, by isotropic, log-concave
measures. It is, perhaps, surprising that extending the classical
gaussian-based results even to natural subgaussian selection
methods, for example, the uniform measure on $\{-1,1\}^n$, is not
simple at all, and in some cases the extension is simply not true.
Moreover, going beyond the subgaussian realm and proving such results for
arbitrary isotropic, log-concave measures is even more difficult, mainly
because the tail estimate that one has for linear functionals is
rather weak. Indeed, in the isotropic, log-concave case the $\psi_1$
and the $\ell_2^n$ norms are equivalent, but
$\|\inr{x,\cdot}\|_{\psi_2}$ might have a strong dependence on the
dimension.

Here, we will study the way a random operator $\Gamma = \sum_{i=1}^N
\inr{X_i,\cdot} e_i$ acts on a convex body, where $(X_i)_{i=1}^N$
are selected according to an isotropic, log-concave measure on
$\R^n$. We will show that many parts of the gaussian theory remain
true for such an operator, with the main difference being that the
classical parameter
$$
\sqrt{n} M^*(K) = \sqrt{n} \int_{S^{n-1}}
\|x\|_{K^\circ} d \sigma \sim \E \sup_{x \in K} \sum_{i=1}^n g_i
x_i
$$
that is used to quantify the phenomena one sees for a
gaussian operator is replaced by $\gamma_2(K,\psi_2)$ (and recall
that $(K,\psi_2)$ is the set of functions $\{\inr{x, \cdot} : x \in
K\}$ endowed with the $\psi_2(\mu)$ norm). Another difference is
that the probabilistic estimates we will obtain for a general
random, isotropic, log-concave operator are much weaker than in the
gaussian or subgaussian cases.

Assume that $K \subset \R^n$ is symmetric. Then $d_{\psi_\alpha}
\sim {\rm diam}(K,\psi_\alpha)$ and $d_{\ell_2^n} \sim {\rm diam}
(K,\ell_2^n)$. For $\alpha=1,2$ and an isotropic measure $\mu$, let
$Q_\alpha(\mu)=\sup_{\theta \in S^{n-1}}
\|\inr{\theta,\cdot}\|_{\psi_\alpha}$ -- the equivalence constant
between the $\psi_\alpha$ norm restricted to linear functionals on
$\R^n$ and the $\ell_2^n$ norm. For example, if $\mu$ is an
isotropic, log-concave measure on $\R^n$ then by Borell's
inequality, $Q_1(\mu) \sim 1$. On the other hand, $Q_2(\mu)$ can
grow polynomially in $n$.

\subsection{The norm of random matrices}
Let $K \subset \R^n$ be a convex body and let $\Gamma : \R^n \to
\R^N$ be the random operator $\sum_{i=1}^N \inr{X_i,\cdot}e_i$,
where $(X_i)_{i=1}^N$ are independent, selected according to an
isotropic, log-concave measure on $\R^n$. Our goal is to estimate
$\E \|\Gamma\|_{K \to \ell_p^N}$, and for the sake of brevity we
will consider the case $p \geq 2$, although the case $1 \leq p <2$
can be handled using similar means.

Let us begin with the relatively simple subgaussian case, when
$Q_2(\mu) \sim 1$.

\begin{Theorem} \label{thm:norm-estimate-subgaussian}
There exists an absolute constant $c$ for which the following holds.
If $p \geq 2$ and $K \subset \R^n$ is a convex body, then for every
integer $N$,
$$
\E \|\Gamma\|_{K \to \ell_p^N} \leq c\left(\gamma_2(K,\psi_2)+
Q_2(\mu) {\rm diam}(K,\ell_2^n) \cdot N^{1/p} \right),
$$
\end{Theorem}

Since the proof of Theorem \ref{thm:norm-estimate-subgaussian} is
rather standard, we will only sketch it here.

\proof Let $p^\prime$ be the conjugate index of $p$. Consider the
random process indexed by $K \times B_{p^\prime}^N$, defined by
$Z_{x,y} = \sum_{i=1}^N \inr{X_i,x}y_i$ and note that for every
$(x,y)$ and $(x^\prime,y^\prime)$,
\begin{equation*}
\|Z_{x,y}-Z_{x^\prime,y^\prime}\|_{\psi_2} \leq d_{\psi_2}
\|y-y^\prime\|_2 + {\rm diam}(B_{p^\prime}^N,\ell_2^N)
\|\inr{X,x-x^\prime}\|_{\psi_2}.
\end{equation*}
Therefore, applying a chaining argument,
$$
\E \sup_{x \in K, \ y \in B_{p^\prime}^N} Z_{x,y} \leq c_1
\left(d_{\psi_2} \gamma_2(B_{p^\prime}^N,\ell_2^N) + {\rm
diam}(B_{p^\prime}^N,\ell_2^N) \gamma_2(K,\psi_2)\right).
$$
To complete the proof, if $G=(g_1,...,g_N)$ is the standard gaussian
vector in $\R^N$ then by the Majorizing Measures Theorem,
$$
\gamma_2(B_{p^\prime}^N,\ell_2^N) \leq c_2 \E \sup_{y \in
B_{p^\prime}^N} \sum_{i=1}^k g_i y_i = c_2 \E \|G\|_{\ell_p^N} \leq
c_3 N^{1/p}.
$$
Also, since $p \geq 2$ then ${\rm diam}(B_{p^\prime}^N,\ell_2^N)=1$
and clearly $d_{\psi_2} \lesssim Q_2(\mu) d_{\ell_2^n}$. Therefore,
$$
\E\|\Gamma\|_{K \to \ell_p^N} \lesssim \gamma_2(K,\psi_2)+Q_2(\mu)
d_{\ell_2^n} N^{1/p},
$$
as claimed.
\endproof

It is simple to verify that Theorem
\ref{thm:norm-estimate-subgaussian} cannot be improved, up to the
constants involved. Indeed, if $\mu$ is the standard gaussian
measure on $\R^n$ then $Q_2(\mu)$ is an absolute constant and
$\gamma_2(K,\psi_2) \sim \gamma_2(K,L_2)=\gamma_2(K,\ell_2^n)$. Let
$(G_i)_{i=1}^N$ be independent copies distributed according to $\mu$
and since $e_1 \in B_{p^\prime}^N$ then
\begin{equation*}
\E \|\Gamma\|_{K \to \ell_p^N} = \E \sup_{x \in K} \sup_{y \in
B_{p^\prime}^N} \sum_{i=1}^N \inr{G_i,x}y_i \geq \E \sup_{x \in K}
\sum_{i=1}^n g_i x_i.
\end{equation*}

Also, if $\|x_0\|_2 = d_{\ell_2^n}$ then
\begin{equation*}
\E \|\Gamma\|_{K \to \ell_p^N} \geq \E\|\Gamma x_0\|_{\ell_p^N} \geq
 c_2 \|x_0\|_2 N^{1/p},
\end{equation*}
showing that the estimate in Theorem
\ref{thm:norm-estimate-subgaussian} is sharp in this case.

Thanks to Theorem B it is possible to replace $Q_2(\mu)$ in Theorem
\ref{thm:norm-estimate-subgaussian} by $Q_1(\mu)$, which, in the
log-concave case, is of the order of an absolute constant.
\begin{Theorem} \label{thm:norm-psi-1-matrix}
There exists an absolute constant $c$ for which the following holds.
Let $K$ be a convex body in $\R^n$. Then for every $p \geq 2$ and
any integer $N$, a random isotropic, log-concave operator $\Gamma$
satisfies that
$$
\E \|\Gamma\|_{K \to \ell_p^N} \leq
c\left(\gamma_2(K,\psi_2) + {\rm diam}(K,\ell_2^n) \cdot
N^{1/p}\right).
$$
\end{Theorem}

\proof Since ${\rm diam}(F,{\psi_1})
\sim {\rm diam}(K,\ell_2^n)$, the claim follows immediately from Theorem B and its extensions to other $\ell_p$ norms for $m=N$
and $F=\{\inr{x,\cdot} : x \in K\}$.
\endproof

An interesting case in which Theorem \ref{thm:norm-psi-1-matrix} can
be used is the ``standard shrinking" phenomenon. Simply
put, standard shrinking is the observation that for every $x \in
\R^n$, and with high probability with respect to the uniform measure
on the Grassman manifold $G_{n,k}$, the random orthogonal projection
$P_E$ satisfies that $\|P_E x\|_2 \leq c \sqrt{k/n} \|x\|_2$. This
property can be extended to a more general situation. Indeed, one
can show that if $K \subset \R^n$ is a convex body, $k^* =
\sqrt{n}M^*(K)/d_{\ell_2^n}$ and $k \geq c_1k^*$, then with high
probability in $G_{n,k}$, ${\rm diam}(P_E K,\ell_2^n) \leq c_2
d_{\ell_2^n} \sqrt{k/n} $. Moreover, this result is sharp, since
Milman's version of Dvoretzky's Theorem (see, for example,
\cite{MilSch}) implies that if $k \leq c_3 k^*$, then with high
probability $P_E K \supset c_4 M^*(K) B_2^k$, and the diameter can
not decrease further.

The shrinking of the diameter for $k \geq k^*$ extends to other
random operators, but even in a relatively simple case, when $\Gamma$
is selected according to the uniform measure on $\{-1,1\}^n$, some
nontrivial machinery is required \cite{Art}, particularly if one
wishes to recover the probabilistic estimate $\sim \exp(-ck)$. The
methods developed in \cite{MPT} (see Corollary 1.9 there) show that
the same is true -- and with the same probability estimate, as long
as $Q_2(\mu) \sim 1$.

Theorem \ref{thm:norm-psi-1-matrix} implies that shrinking does
happen for a random isotropic, log-concave operator -- though with
a weaker probabilistic estimate. Indeed, consider the operator
$A=\Gamma/\sqrt{n}$, let $K \subset \R^n$ be a convex body and set
$k^\prime = \gamma_2(K,\psi_2)/d_{\ell_2^n}$. Then, with high
probability,
$$
{\rm diam}(AK,\ell_2^k) \lesssim \frac{1}{\sqrt{n}} \left(\gamma_2(K,\psi_2)+{\rm
diam}(K,\ell_2^n) \sqrt{k} \right) \lesssim \sqrt{\frac{k}{n}} {\rm
diam}(K,\ell_2^n),
$$
as long as $k \geq k^\prime$. Since
$$
 c_1\sqrt{n} M^*(K) \leq
\gamma_2(K,\psi_2) \leq c_2 Q_2(\mu) \sqrt{n} M^*(K),
$$
it follows that if $\mu$ happens to be subgaussian, i.e. if $Q_2(\mu) \sim 1$,
then $k^\prime$ and $k^*$ are equivalent.

\subsection{Low--$M^*$ estimates}
Given a convex body $K \subset \R^n$ and $k \leq n$, one would like
to find a subspace $E \subset \R^n$ for which the Euclidean
diameter of $K \cap E$ is as small as possible. We refer the reader
to \cite{Mil,MPT} for a brief description of the progress made on
this problem.

In \cite{PT,PT1} it was shown that if $E$ is the kernel of a random
orthogonal projection (or of a gaussian projection), and if
\begin{equation} \label{eq:r_k}
r_N^* =\inf \left\{ r>0 : \sqrt{n}M^*(K \cap rS^{n-1})/{\sqrt{N}}
\leq {cr} \right\},
\end{equation}
then ${\rm diam}(E \cap K) \leq r^*_N$, where $c$ is an absolute
constant.

Since the original proof of this result is based on the structure of
gaussian variables or that of the Haar measure on $G_{n,k}$, extending it
to other natural random operators is not trivial. Equation
\eqref{eq:r_k} was extended to the subgaussian case in \cite{MPT}
using a subgaussian version of Theorem A. It was shown that if $\mu$
is isotropic and $\Gamma=\sum_{i=1}^N \inr{X_i, \cdot} e_i$ (with
$X_1,...,X_N$, independent, distributed according to $\mu$), then
with high probability,
\begin{equation} \label{eq:r_k-1}
 {\rm diam}(K \cap {\rm ker} \Gamma)
\leq \inf\left\{ r >0 : Q_2(\mu) \gamma_2(K \cap r
S^{n-1},\psi_2)/{\sqrt{N}} \leq cr \right\}.
\end{equation}
Therefore, if $\mu$ is isotropic and $Q_2(\mu) \sim 1$, (that is, if $\Gamma$ is an isotropic, subgaussian operator) then \eqref{eq:r_k} is true. Applying Theorem
A, the fact that for an isotropic, log-concave measure $Q_1(\mu) \sim 1$ and the proof from \cite{MPT}, one has
\begin{Theorem} \label{thm:small-diameter-applications}
There exist absolute constants $c$ and $c_1$ for which the following
holds. Let $\Gamma:\ell_2^n \to \ell_2^N$ be a random isotropic
log-concave operator. Then for a convex body $K \subset \R^n$ one
has
\begin{equation*}
\E \left({\rm diam}(K \cap {\rm ker} \Gamma) \right) \leq c_1
\inf\left\{ r
>0 :  \gamma_2(K \cap r S^{n-1},\psi_2)/{\sqrt{N}} \leq cr
\right\},
\end{equation*}
and a similar estimate holds with high probability.
\end{Theorem}

Again, Theorem \ref{thm:small-diameter-applications} extends the
classical result to any isotropic log-concave case ensemble, with $\gamma_2(K,\psi_2)$
taking the place of $\sqrt{n}M^*(K)$ -- though with a weaker
probabilistic estimate.

\subsection{The process indexed by $S^{n-1}$} \label{sec:sampling}
This section is devoted to a problem that is far from being fully
solved -- the behavior of the process
\begin{equation} \label{eq:process-sphere}
\sup_{\theta \in S^{n-1}} \left|\frac{1}{N}\sum_{i=1}^N
\inr{X_i,\theta}^2-1\right|,
\end{equation}
where $X_1,...,X_N$ are selected independently according to an
isotropic, log-concave measure on $\R^n$.

In \cite{ALPT} the authors solved the following facet of this
problem: Given $\eps>0$ and $0<\delta<1$, how many random points
$X_1,...,X_N$ are needed to ensure that with probability $1-\delta$,
$$
\sup_{\theta \in S^{n-1}} \left|\frac{1}{N}\sum_{i=1}^N
\inr{\theta,X_i}^2 -1 \right| < \eps?
$$
An equivalent formulation of this question is to find the smallest
$N$ that would still guarantee that a random, isotropic, log-concave
operator $\Gamma$ embeds $\ell_2^n$ in $\ell_2^N$ $1+\eps$
isomorphically.

This problem has been studied extensively in recent years (e.g.
\cite{KLS,Bour,GiaMil,Rud,GueRud,Paou,Men:weak,Aub}), in which the
estimate has been improved from the initial $N \geq
c(\eps,\delta)n^2$ in \cite{KLS} to the best possible estimate of $N
\geq c(\eps,\delta)n$, proved in \cite{ALPT}. In fact, what was
actually proved in \cite{ALPT} is the following:

\begin{Theorem} \label{thm:ALPT}
There exist absolute constants $C$, $c$ and $c_1$ for which the following
holds. Let $\mu$ be an isotropic, log-concave measure on $\R^n$
and let $(X_i)_{i=1}^N$ be independent, distributed according to
$\mu$. Then, for every $t \geq 1$ and every $1 \leq N \leq
\exp(\sqrt{n})$, with probability at least $1-2\exp(-ct\sqrt{n})$,
for every $I \subset \{1,...,N\}$,
\begin{equation} \label{eq:ALPT1}
\sup_{\theta \in S^{n-1}} \left(\sum_{i \in I} \inr{\theta,X_i}^2
\right)^{1/2} \leq C\left(\sqrt{n} + \sqrt{|I|}\log(eN/|I|) \right).
\end{equation}
Moreover, for every $c_1n \leq N \leq \exp(\sqrt{n})$ and every $s,t
\geq 2$,
\begin{equation} \label{eq:ALPT2}
\sup_{\theta \in S^{n-1}} \left|\frac{1}{N}\sum_{i=1}^N
\inr{\theta,X_i}^2 -1 \right| \leq
C\left(ts\sqrt{\frac{n}{N}}\log(eN/n) + s \frac{n}{N} \right)
\end{equation}
with probability at least $1-2\exp(-cs\sqrt{n})-2\exp(-c
\min\{u,v\})$, where $u=t^2s^2n\log^2(eN/n)$ and
$v=(t/s)\sqrt{nN}/\log(eN/n)$.
\end{Theorem}

Although Theorem \ref{thm:ALPT} beautifully resolves the case $N
\sim n$, its proof has certain weaknesses from the point of view of
empirical processes theory and the general understanding of the
process \eqref{eq:process-sphere}. First of all, \eqref{eq:ALPT2} is
derived from \eqref{eq:ALPT1} using a decomposition and contraction
argument, just like our Theorem C is derived from Theorem B. Hence,
there is an intrinsic logarithmic looseness in \eqref{eq:ALPT2} -- a
superfluous factor of $\log N$ for $N \geq c(\beta)n^{1+\beta}$ for
any $\beta>0$.

Second, the proof of Theorem \ref{thm:ALPT} relies on the Euclidean
nature of the problem in a very strong way: that the given class is
a class of linear functionals on $\R^n$, that the indexing set is
the entire sphere and that the measure is isotropic, log-concave (in
particular, that $Q_1(\mu) \sim 1$ and that the Euclidean norm of a
random point concentrates around $\sqrt{n}$). Hence, the method of
\cite{ALPT} cannot be extended beyond this limited setup, even to
obtain an analogous result for a small subset of the sphere as an
indexing class. Naturally, it is also impossible to obtain an
``empirical processes" result like Theorem A in this way. A
consequence of this limitation is that the method of \cite{ALPT}
cannot be used to prove the applications presented in the two previous
sections (i.e., estimates on the norm $\|\Gamma\|_{K \to \ell_2^N}$,
the shrinking phenomenon, and low-$M^*$ estimates) since those
applications require accurate information on the way $\Gamma$ acts
on arbitrary subsets of $\R^n$ rather than on the entire sphere.

Process \eqref{eq:process-sphere} is very far from being understood
when one goes beyond the case $N \sim n$. A reasonable conjecture is
that for any $N \gtrsim n$, with high probability/in expectation,
\begin{equation} \label{eq:sampling}
\sup_{\theta \in S^{n-1}} \left|\frac{1}{N} \sum_{i=1}^N
\inr{\theta,X_i}^2 - 1 \right| \leq c\sqrt{\frac{n}{N}},
\end{equation}
which is the situation for the gaussian ensemble.

Below, we will indicate some of the problems one faces when trying
to verify this conjecture, with the main one being that very little
is known on the metric structure endowed on $S^{n-1}$ by a
log-concave measure.

Currently, the best estimate on \eqref{eq:process-sphere} in the
range $c(\beta)n^{1+\beta} \leq N \leq \exp(\sqrt{n})$ for any $\beta
>0$ is $c_1(\beta)\sqrt{(n\log n)/N}$. This is a corollary of Theorem A,
and the {\it suboptimal} estimate from \cite{Men:weak}, that $\gamma_2(S^{n-1},\psi_2) \lesssim
\sqrt{n \log n}$ for $\mu$ that is supported in $c_2\sqrt{n}B_2^n$
(the so-called {\it small diameter case}). Note that the small diameter assumption can be made without loss of generality as long as $N \leq \exp(\sqrt{n})$ thanks to the result of Paouris \cite{Paou} which
states that for $N \leq \exp(\sqrt{n})$, $\E \max_{i \leq N}
\|X_i\|_2 \lesssim \sqrt{n}$. Hence, for those values of $N$, one
may assume that $\mu$ is supported in $c_2\sqrt{n}B_2^n$, implying that if $c(\beta)n^{1+\beta} \leq N \leq \exp(\sqrt{n})$ then Theorem A improves Theorem \ref{thm:ALPT} and gives the best known estimate on \eqref{eq:process-sphere}.

We believe that under the small diameter assumption, the extra
logarithmic term in $\gamma_2(S^{n-1},\psi_2)$ could be removed.
Indeed, if $\mu$ is supported on a ball of radius $\sim\sqrt{n}$,
``most" directions $\theta \in S^{n-1}$ have a $\psi_2$ norm that is
bounded by an absolute constant (see, for example,
\cite{Gia:survey}). Unfortunately, even under a small diameter
assumption, there is very little information on the geometry of the
set of these ``good" directions, except that it is a very large
subset of the sphere.

The second step towards a complete solution, and most likely the
more difficult one, is when $N \geq \exp(\sqrt{n})$. Here, one can
no longer assume that $\mu$ is supported in a ball of radius $\sim
\sqrt{n}$, and thus both Theorem \ref{thm:ALPT} and the bound on
$\gamma_2(S^{n-1},\psi_2)$ from \cite{Men:weak} fail. Moreover, when leaving the small diameter case, it
is not known whether there is even a {\it single} direction $\theta$
for which $\|\theta\|_{\psi_2} \sim 1$.

\subsubsection{The unconditional case} We end this note with an
example of how the $\sqrt{\log{n}}$ factor may be removed in a
special case, when $\mu$ is unconditional. This example illustrates
the difficulties that one is likely to encounter in the general
case, where there is little structure at our disposal.

The argument has two parts. First, we will show that one may
consider a slightly different ``small diameter" assumption, and
second, that under this assumption, the metric entropy
$\log N(S^{n-1},\eps B_{\psi_2})$ is well behaved.

For the first part, note that by the Bobkov-Nazarov Theorem
\cite{Bob-Naz,Gia:survey}, if $N \sim n^\alpha$ and if we denote the
$j-th$ coordinate of a monotone rearrangement of the coordinates of the vector $X_i$
by $(X_i)_j^*$,
then with high probability,
for every $1 \leq i \leq N$ and $1 \leq j \leq n$, $(X_i)_j^* \leq
c_\alpha \log(en/j)$. Hence, without loss of generality we may
assume that $\mu$ is supported in $c_1(\alpha) B_{\psi_1^n}$. This gives more accurate information than the standard small diameter assumption, that $\mu$ is supported in $c\sqrt{n}B_2^n$. In particular, we may assume that almost surely, for every $j \leq n$, $(X)^*_j \leq c_\alpha \log(en/j)$. Since $\mu$ is unconditional, then for every $\theta \in S^{n-1}$ the random variable $\inr{X,\theta}$ has
the same distribution as $\sum_{j=1}^n \eps_j |\inr{X,e_j}|\theta_j$, where $(\eps_j)_{j=1}^n$ are i.i.d. Bernoulli random variables. Hence, for any $p \geq 1$,
\begin{align*}
& \left(\E_X |\inr{X,\theta}|^p \right)^{1/p} =  \left(\E_{X \times \eps}
\left|\sum_{j=1}^n \eps_j |\inr{X,e_j}| \theta_j \right|^p \right)^{1/p}
\\
\leq & c\sqrt{p} \left(\E_X \left|\sum_{j=1}^n |\inr{X,e_j}|^2 \theta_j^2
\right|^{p/2} \right)^{1/p}
\leq c\sqrt{p}\left(\E_X \left|\sum_{i=1}^n (X^2)_j^* (\theta^2)_j^* \right|^{p/2} \right)^{1/p}
\\
\leq & c_1(\alpha) \sqrt{p} \left(\sum_{j=1}^n (\theta^2)_j^* \log^2(en/j)
\right)^{1/2}.
\end{align*}

In particular, for every $\theta \in S^{n-1}$,
$$
\|\theta\|_{\psi_2}
\leq c(\alpha) \left(\sum_{j=1}^n (\theta^2)_j^* \log^2(en/j)
\right)^{1/2}
$$
and ${\rm diam}(S^{n-1},\psi_2) \leq c(\alpha) \log
n$.

Now, just as in \cite{Men:weak} one may show that for every $\eps
\leq 2$, the covering numbers satisfy $N(S^{n-1}, \eps B_{\psi_2})
\leq (c_2/\eps)^n$. Thus, it remain to estimate the covering numbers
for larger scales.

To that end, we will use a minor modification of the sets $N_\ell$
and $B_m$ that appeared in Section \ref{sec:diam}.

Let $ A_\ell = \left\{z \in B_2^n : \ |{\rm supp} (z)| \leq \ell, \
\|z\|_\infty \leq 1/\sqrt{\ell} \right\}$, fix $r$ such that $2^r
\leq n/10$ and let $\eps_r = \log(en/2^r)$. Set $N_{2^j} \subset
A_{2^j}$ to be an $\eps_r (2^j/n)$-cover of $A_{2^j}$ with respect
to the $\ell_2^n$ norm and define
\begin{equation*}
B_r = \left\{ z \in B_2^n : \ |{\rm supp}(z)| \leq 2^r, \ {\rm
supp}(z)=\bigcup_{j=0}^{r-1} I_j \ , \  P_{I_j} z \in N_{2^j}
\right\},
\end{equation*}
where $I_j$ are disjoint sets of coordinates with $|I_0|=2$ and
$|I_j|=2^j$ for $j \geq 1$.

It is standard to verify that $|B_r| \leq \exp(c_0 2^r \log
(en/2^r))$ and that for every $\theta \in S^{n-1}$ there is some
$\tilde{\theta} \in B_r$ whose support is denoted by $I$, such that
\begin{align*}
\|\theta - \tilde{\theta}\|_{\psi_2} \leq & c_1\left(
\sum_{j=0}^{r-1} \|P_{I_j} (\theta - \tilde{\theta})\|_2
\log(en/2^j) + \|P_{I^c}\theta\|_2 \log(en/2^r) \right)
\\
\leq & c_1 \left( \frac{\eps_r}{n} \sum_{j=0}^{r-1} 2^j \log(en/2^j)
+ \log(en/2^r) \right) \leq c_2 \eps_r
\end{align*}
for $c_1$ and $c_2$ that depend on $\alpha$.

Therefore, $B_r$ is a $c_2\eps_r$-cover of $S^{n-1}$ with respect to
the $\psi_2$ norm, implying that
$$
\log N(S^{n-1}, c_2\eps_r B_{\psi_2}) \leq c_0 2^r \log (en/2^r).
$$
It is well known \cite{Tal:book} that if $(T,d)$ is a metric space
then
$$
\gamma_2(T,d) \lesssim \int_0^{{\rm diam}(T,d)}
\sqrt{\log(N(T,\eps,d))}d\eps,
$$
and thus a simple calculation of this entropy integral shows that
$$\gamma_2(S^{n-1},\psi_2) \leq c_3(\alpha)
\sqrt{n},
$$
proving our claim.
\endproof

\newpage
\footnotesize {

\end{document}
\begin{thebibliography}{10} \frenchspacing
%
\bibitem{ALPT} R. Adamczak, A. Litvak, A. Pajor, N.
Tomczak-Jaegermann, Quantitative estimates of the convergence of the
empirical covariance matrix in log-concave ensembles, J. Amer. Math. Soc. 23 535-561, 2010.

\bibitem{ALPT2} R. Adamczak, A. Litvak, A. Pajor, N.
Tomczak-Jaegermann,Restricted Isometry Property of Matrices with
Independent Columns and Neighborly Polytopes by Random Sampling,
preprint.

\bibitem{Art} S. Artstein, Change in the diameter of a convex body
under a random sign-projection, Geometric aspects of functional
analysis, Lecture Notes in Math. 1850, 31--39, Springer, 2004.

\bibitem{Aub} G. Aubrun, Sampling convex bodies: a random matrix
approach, Proc. Amer. Math. Soc. 135, 1293-1303, 2007.

\bibitem{BGMN} F. Barthe, O. Gu\'{e}don, S. Mendelson, A. Naor,
A probabilistic approach to the geometry of the $\ell_p^n$ ball,
Ann. Probab. 33(2), 480-513, 2005.

\bibitem{BMN} P.L. Bartlett, S. Mendelson, J. Neeman,
$\ell_1$-regularized linear regression: Persistence and oracle
inequalities, submitted.

\bibitem{Bob-Naz} S.G. Bobkov, F.L. Nazarov, On convex bodies and
log-concave probability measures with unconditional basis, Geometric
Aspects of Functional Analysis, Lecture Notes in Mathematics 1807,
53-69, 2003.

\bibitem{Bor} C. Borell, The Brunn-Minkowski inequality in Gauss
space, Invent. Math. 30 207-216, 1975.

\bibitem{Bour} J. Bourgain, Random points in isotropic convex
bodies, in Convex Geometric Analysis (Berkeley, CA, 1996) Math. Sci.
Res. Inst. Publ. {34} (1999), 53-58.

\bibitem{CRT} E. Candes, J. Romberg, T. Tao, Robust uncertainty
principles: exact signal reconstruction from highly incomplete
frequency information, IEEE Trans. Info. Theory 52 (2), 489--509,
2006.

\bibitem{CT1} E. Candes, T. Tao, Near-optimal signal
recovery from random projections: universal encoding strategies?
IEEE Trans. Inform. Theory 52 (12), 5406--5425, 2006.

\bibitem{CT2} E. Candes, T. Tao, The Dantzig selector:
statistical estimation when $p$ is much larger than $n$, Ann.
Stat. 35 (6), 2313--2351, 2007.

\bibitem{Dud-book} R. M. Dudley, {\it Uniform Centra Limit
Theorems}, Cambridge Studies in Advanced Mathematics 63, Cambridge
University Press, 1999.

\bibitem{F} X. Fernique, R\'{e}gularit\'{e} des trajectoires des fonctiones al\'{e}atoires
gaussiennes, Ecole d'Et\'{e} de Probabilit\'{e}s de St-Flour 1974,
Lecture Notes in Mathematics 480, 1-96, Springer-Verlag 1975.


\bibitem{GiaMil} A.A. Giannopoulos, V.D. Milman, Concentration
property on probability spaces, Adv. Math. 156, 77-106, 2000.

\bibitem{Gia:survey} A.A. Giannopoulos, Notes on isotropic convex
bodies, available at http://users.uoa.gr/$\sim$apgiannop/

\bibitem{GZ84} E. Gin\'{e} and J. Zinn, Some limit theorems for empirical
processes, Ann. Probab. 12(4), 929-989, 1984.

\bibitem{GidlP} E. Gin\'{e}, V. H. de la Pe\~{n}a, {\it Decoupling},
Springer, 1999.

\bibitem{GR} E. Greenshtein, Y. Ritov, Persistence in high-dimensional
linear predictor selection and the virtue of overparametrization,
Bernoulli, 10(6), 971--988, 2004.


\bibitem{GueRud} O. Gu\'{e}don, M. Rudelson, $L_p$ moments of
random vectors via majorizing measures, Adv. Math. 208(2), 798-823,
2007.

\bibitem{KLS} R. Kannan, L. Lov\'{a}sz, M. Simonovits, Random
walks and $O^*(n^5)$ volume algorithm for convex bodies, Random
structures and algorithms, 2(1) 1-50, 1997.

\bibitem{KM} B. Klartag, S. Mendelson,
Empirical Processes and Random Projections, Journal of Functional
Analysis, 225(1) 229-245, 2005.


\bibitem{LT} M. Ledoux, M. Talagrand,  {\it Probability in
  Banach spaces. Isoperimetry and processes}, Ergebnisse der Mathematik
  und ihrer Grenzgebiete (3), vol. 23.  Springer-Verlag, Berlin, 1991.

\bibitem{MPT} S. Mendelson, A. Pajor, N.
Tomczak-Jaegermann, Reconstruction and subgaussian operators,
Geometric and Functional Analysis, 17(4), 1248-1282, 2007.

\bibitem{Men:weak} S. Mendelson, On weakly bounded empirical processes,
Math. Annalen, 340(2), 293-314, 2008.

\bibitem{Mil} V.D. Milman, A note on a low-$M^*$ estimate, in {\it
Geometry of Banach spaces (strobl, 1989)}, London Math. Soc. Lecture
Notes Ser. 158, 219-229, Cambridge University Press, 1990.

\bibitem{MilSch} V.D. Milman, G. Schechtman, { Asymptotic theory of
finite dimensional normed spaces}, Lecture Notes in Mathematics
1200, Springer, 1986.
%
\bibitem{PT} {A. Pajor, N. Tomczak-Jaegermann,}  Subspaces
  of small codimension of finite-dimensional Banach spaces,
  Proceedings of the AMS. 97(4),
  637-642, 1986.

\bibitem{PT1} {A. Pajor, N. Tomczak-Jaegermann,}  Nombres
  de Gelfand et sections euclidiennes de grande dimension.
  (French) [Gelfand numbers and high-dimensional Euclidean sections]
  S\'{e}minaire d'Analyse Fonctionelle 1984/1985, Publ. Math.
  Univ. Paris VII, 26, Univ. Paris VII, Paris,  37--47, 1986.

\bibitem{Paou} G. Paouris, Concentration of mass on convex bodies,
Geometric and Functional Analysis, 16(5), 1021-1049, 2006.


\bibitem{Pis:book} G. Pisier,
{\it The volume of convex bodies and Banach space geometry},
Cambridge University Press, 1989.

\bibitem{Rud} M. Rudelson, Random vectors in the isotropic
position, J. Funct. Anal. 164, 60-72, 1999.

\bibitem{Rud-selectors} M. Rudelson, Extremal distances between
sections of convex bodies, Geom. Funct. Anal.
14(5), 1063-1088, 2004.

\bibitem{Tal87} M. Talagrand, Regularity of Gaussian processes,
Acta Math. 159, 99-149, 1987.


\bibitem{Tal:AJM94} M. Talagrand, The supremum of some canonical
processes, American Journal of Mathematics 116, 283--325, 1994.

\bibitem{Tal:book} {M. Talagrand,} {\it The generic chaining},
Springer, 2005.

\bibitem{VW}{A.W. Van der Vaart, J.A. Wellner, }
{\it Weak convergence and empirical processes}, Springer Verlag,
1996.
%

\end{thebibliography}
